\documentclass[12pt]{article}
\usepackage{amsmath,amssymb,amsfonts,amscd,amsxtra}
\usepackage{latexsym}

\input{epsf.sty}
\usepackage{epsfig}
\usepackage{graphicx,color}
\usepackage{graphics}
\usepackage{subfigure}
\usepackage{appendix}
\usepackage[top=1in, bottom=1in, left=1in, right=1in]{geometry}

\renewcommand{\theequation}{\thesection.\arabic{equation}}

\newtheorem{thm}{Theorem}[section]

\newtheorem{lem}[thm]{Lemma}
\newtheorem{prop}[thm]{Proposition}
\newtheorem{rmk}[thm]{Remark}

\renewcommand{\Im}{{\mbox{Im}}}
\renewcommand{\Re}{{\mbox{Re}}}

\newcommand{\norm}[1]{\left\Vert#1\right\Vert}
\newcommand{\abs}[1]{\left\vert#1\right\vert}

\newcommand{\qed}{\hfill \ensuremath{\square}}

\renewcommand\appendix{\par
  \setcounter{section}{0}
  \setcounter{subsection}{0}
  \setcounter{figure}{0}
  \setcounter{table}{0}
  \renewcommand\thesection{Appendix \Alph{section}}
  \renewcommand\theequation{\Alph{section}.\arabic{equation}}
  \renewcommand\thefigure{\Alph{section}.\arabic{figure}}
  \renewcommand\thetable{\Alph{section}.\arabic{table}}
  \renewcommand\thethm{\Alph{section}.\arabic{thm}}
}

\numberwithin{equation}{section}

\date{}

\title{Fano resonance in metallic grating via strongly coupled subwavelength resonators}
\author{
Junshan Lin\thanks{\footnotesize Department of Mathematics and Statistics, Auburn University, Auburn, AL 36849 (jzl0097@
auburn.edu). Junshan Lin was partially supported by the NSF grant DMS-1719851.}
 \; and Hai Zhang\thanks{\footnotesize 
 Department of Mathematics, 
  HKUST,  Clear Water Bay, Kowloon, Hong Kong SAR, China (haizhang@ust.hk). Hai Zhang was supported by Hong Kong RGC grant GRF 16304517 and GRF 16306318.}}

\begin{document}

\maketitle

\begin{abstract}
We investigate the Fano resonance in grating structures by using coupled resonators.
The grating consists of a perfectly conducting slab with periodically arranged subwavelength slit holes, 
where inside each period, a pair of slits sit very close to each other.
The slit holes act as resonators and are strongly coupled. 
It is shown rigorously that there exist two groups of resonances corresponding to poles of the scattering problem.
One sequence of resonances have imaginary part on the order of $\varepsilon$, where $\varepsilon$ is the size of the slit aperture,  while the other sequence have imaginary part on the order of $\varepsilon^2$.
When coupled with the incident wave at resonant frequencies, 
the narrow-band resonant scattering induced by the latter will interfere with the broader background resonant radiation induced by the former. 
The interference of these two resonances generates the Fano type transmission anomaly,
which persists in the whole radiation continuum of the grating structure as long as the slit aperture size is small compared to the incident wavelength.
\end{abstract}


\setcounter{equation}{0}
\setlength{\arraycolsep}{0.25em}
\section{Introduction}\label{sec:introduction}
Fano resonance is a special type of resonant wave phenomenon that gives rise to an asymmetric spectral line shape as opposed to the conventional symmetric one.
Typically the corresponding transmission or reflection signal exhibits a sharp transition from peak to dip within a narrow band,
and this interesting phenomenon has been explored extensively in a variety of oscillating systems and wave structures since first discovered by Ugo Fano in \cite{fano}.
In general, Fano resonance is induced by the constructive or destructive interference between a broad background radiation and a narrow resonant state.
This can be realized when a discrete localized state becomes coupled to a continuum of states, as in the atomic system investigated in the original paper by Ugo Fano 
 \cite{fano} and many other optical structures including dielectric spheres, photonic crystals, plasmonic metasurfaces, etc \cite{hsu3, limonov, lukyanchuk}.
Another way to realize Fano resonance is by using two coupled resonators with distinct resonance strength so that
the stronger resonant state would intefere with the weaker one.
See, for instance, \cite{fan1, li, limonov, TKT, zhou} for the realization of Fano resonance using coupled oscillators, cavities, micro-spheres, micro-rings, etc.

For photonic grating structures comprising an array of holes etched in a dielectric or metallic slab, it is known that
certain bound states or discrete localized states appear in the continuum (sometimes called BIC), which can lead to Fano resonance.
In more details, there exist eigenvalues embedded in the continuous (radiation) spectrum of the underlying differential operators,
and the bound states associated with those eigenvalues become narrow-band resonant states that 
interact with the background radiation when the symmetry of the system is destroyed.
Such interaction generates the asymmetric resonance transmission line shape, and the principle is the same as the one in  \cite{fano} and others mentioned above.
We refer the reader to \cite{bugakov,fan2,hsu2,hsu3,shipman05,shipman10,lu2} and the references therein for the theoretical and experimental investigation of BICs
in grating structures.

BICs have been rigorously proven to exist in the first diffraction continuum for grating structures where only one diffraction mode propagates \cite{bonnet_starling94,shipman10}. 
Such BICs are usually preserved by the symmetry in the geometry (such as reflection or rotational invariant of the grating) and the incident wave pattern (normal incidence).
Hence Fano resonance naturally appears in this particular diffraction region under a small perturbation of the system which destroys the symmetry.
The related mathematical studies was first carried out in \cite{shipman05,shipman10}, where boundary integral equations and perturbation theories
are applied to explain and predict the resonant transmission peaks and dips.
In connection to this, in \cite{lin_shipman_zhang} we provided quantitative analysis for embedded eigenvalues 
and their perturbations as resonances in the context of periodic metallic grating structures with small holes,
see more discussions in Section 1.2.
Outside the first diffraction regime, however, it is very challenging to find the BICs and the associated Fano resonances, see 
\cite{bugakov,porter,hsu2,lu1,lu2} for the recent attempts in the direction.

In this paper, instead of relying on BICs for Fano resonance, 
we employ the second approach mentioned above and investigate mathematically the realization of Fano resonance in grating structures by using coupled resonators.
Although such an idea has been explored empirically for a variety of coupled resonators \cite{fan1, li, limonov, TKT, zhou}, a rigorous mathematical analysis has not been given yet. It is the goal of this paper to fill the gap here. 
We shall consider a grating structure which consists of a pair of two strongly coupled slit hole resonators in each period.
In contrast to usual grating structures where Fano resonance occurs above the first radiation continuum only for very special configurations, 
the proposed photonic structures allow for the appearance of Fano resonance in the whole radiation continuum,  as long as the size of the opening
for the slit hole resonators $\varepsilon$ is small compared to the incident wavelength $\lambda$.
More precisely, we shall prove rigorously that the coupled subwavelength resonant slit holes would generate two group of resonances with close resonant frequencies
(real parts of resonances) but distinct resonance strength (imaginary parts of resonances).
The first group of resonances attain an imaginary part on an order of $\varepsilon$, while the second attain much smaller imaginary parts on an order of $\varepsilon^2$. 
Consequently, when coupled with the external radiation at resonant frequencies, 
the narrow-band resonant scattering induced by the latter will interfere with the broader resonant scattering induced by the former.
This gives rise to the Fano type transmission anomaly near the resonant frequencies.
We point out a closely related work in \cite{eric10-2}, where it is shown that  two group of resonances induced by a non-periodic structure with two coupled cavity holes would attain
 imaginary parts on the order of $\varepsilon$ and $\varepsilon^3$ respectively. 
We also refer the reader to \cite{habib2019} for using strongly coupled bubble dimers to generate double-negative refractive index media.

\subsection{Problem setup}
\begin{figure}
\begin{center}
\includegraphics[height=6cm,width=14cm]{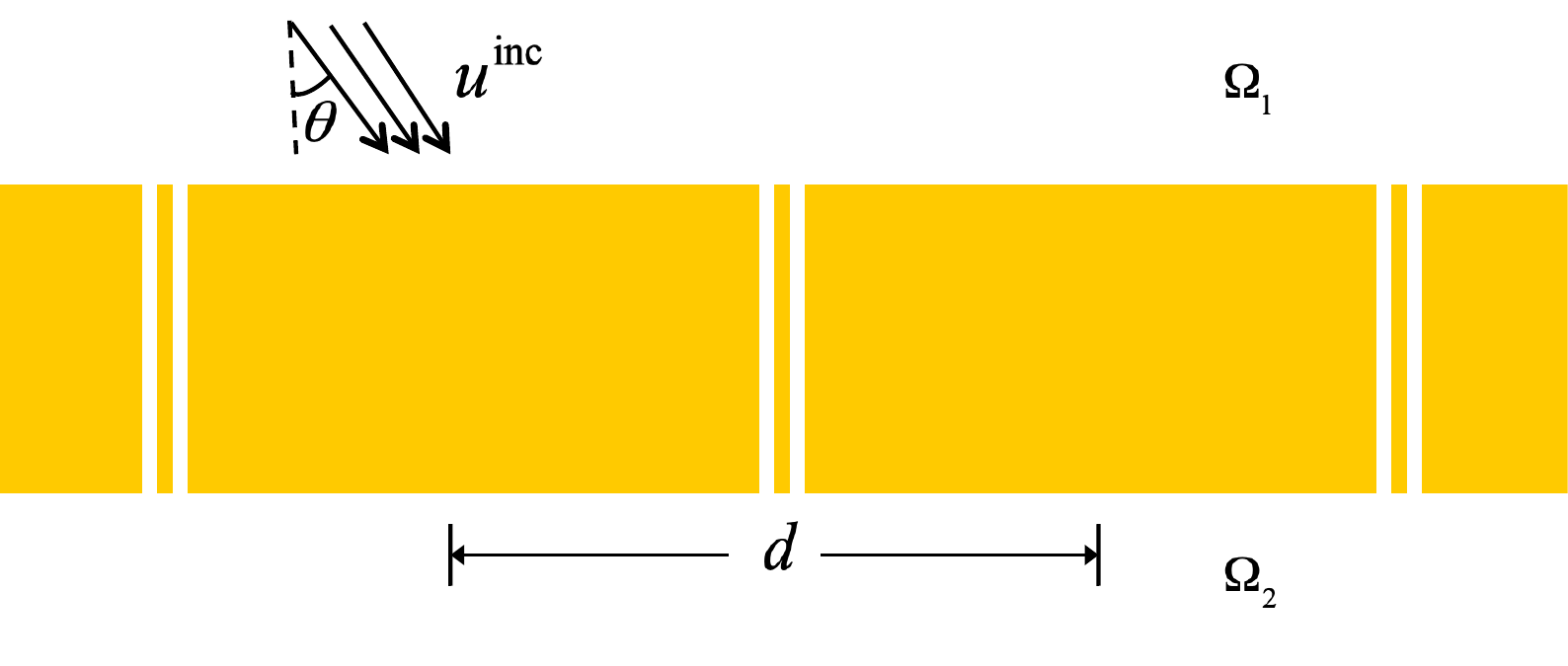}
\vspace*{-20pt}
\caption{Metallic grating structure where each period consists of two narrow slits $S_\varepsilon^{0,-}$ and $S_\varepsilon^{0,+}$ close to each other.
Each slit has length $1$ and width $\varepsilon$.  The upper and lower aperture of the slit $S_{\varepsilon}^{0,\pm}$ is denoted as $\Gamma^{\pm}_{1,\varepsilon}$ and $\Gamma^{\pm}_{2,\varepsilon}$, respectively. The domain exterior to the perfect conductor is denoted as $\Omega_{\varepsilon}$, which consists of the slit region $S_\varepsilon$, the domains above the slab $\Omega_1$, and the domain below the slab $\Omega_2$.
 }\label{fig:prob_geo}
\end{center}
\end{figure}

The metallic grating structure investigated in this paper is shown in Figure~\ref{fig:prob_geo}, where each period contains
a pair of identical narrow slits close to each other. 
Throughout the paper, we denote the grating period by $d$ and the reciprocal lattice constant by $b=2\pi/d$. It is assumed that $d>1$ so that no real-valued eigenvalues would exist below the continuous spectrum 
of the periodic differential operator (cf. Proposition 4.6 and Remark 4.2 in \cite{lin_zhang18_1}). This assumption allows us to focus on resonances lying in the continuous spectrum only but it is not essential.
The slits are invariant along the $x_3$ direction,
and they occupy the region $\displaystyle{S_\varepsilon=\bigcup_{n=0}^{\infty} (S_\varepsilon^{0,-} \cup S_\varepsilon^{0,+}+ nd)}$ in the $x_1x_2$ plane, where the two subwavelength slits are given by
$$\displaystyle S_\varepsilon^{0,\pm}:=\left\{(x_1,x_2)\;|\; (\pm \ell - 1) \frac{ \varepsilon}{2}  <x_1< (\pm \ell + 1)  \frac{ \varepsilon}{2} ,\; 0<x_2<1 \right\}. $$
In the above, $\ell=O(1)$ is some positive constant greater than $1$ such that the distance between two slits is given by $(\ell-1)\varepsilon$.
Denote the upper and lower apertures of the slit $S_{\varepsilon}^{0,\pm}$ by $\Gamma^{\pm}_{1,\varepsilon}$ and $\Gamma^{\pm}_{2,\varepsilon}$.
Let $\Omega_1$ and $\Omega_2$ be the domain above and below the grating respectively, and $\Omega_\varepsilon$ be
the exterior domain of the metallic structure.

We consider the time-harmonic transverse magnetic (TM) situation.  The incident magnetic field is perpendicular to the  $x_1x_2$ plane and its $x_3$ component is the scalar function
\begin{equation}
  u^\mathrm{inc}(x) = e^{i k( x_1\sin\theta \,-\, (x_2-1)\cos \theta )} = e^{i\kappa x_1-i\zeta_0 (x_2-1)},
\end{equation}
in which $k$ is the free-space wavenumber, $\theta\in (-\pi/2,\pi/2)$ is the angle of incidence, $\kappa= k \sin \theta$ is the Bloch wavenumber, and $\zeta_0=\sqrt{k^2-\kappa^2}>0$.  The total field $u_\varepsilon(x)$ satisfies the Helmholtz equation
\begin{equation}\label{eq:Helmholtz}
\Delta u_\varepsilon + k^2 u_\varepsilon \;=\; 0 \quad \mbox{in} \; \Omega_\varepsilon,
\end{equation}
and the Neumann boundary condition
\begin{equation}\label{eq:Neumann}
\frac{\partial u_{\varepsilon}}{\partial \nu} = 0  \quad \mbox{on} \; \partial \Omega_{\varepsilon},
\end{equation}
where $\nu$ is the unit normal vector pointing to $\Omega_\varepsilon$. In addition, the quasi-periodicity of the solution is consistent with the incident field such that
\begin{equation}\label{eq:quasi-periodic}
u_\varepsilon(x_1+d,\,x_2) \;=\; e^{i\kappa d}u_\varepsilon(x_1,x_2).
\end{equation}
Finally, the diffracted field is outgoing and $u_\varepsilon(x)$ admits the Rayleigh-Bloch expansions (cf.~\cite{bao95, bonnet_starling94, shipman10}) above and below the grating structure:
\begin{eqnarray}\label{eq:rad_cond}
  && u_\varepsilon(x_1,x_2) \;=\; u^\mathrm{inc}(x_1,x_2) + \sum_{n=-\infty}^{\infty} u_{n,1}^{s} e^{i \kappa_n x_1 + i\zeta_n  x_2 } \quad \mbox{in} \; \Omega_1, \label{eq:rad_cond1}\\
  && u_\varepsilon(x_1,x_2) \;=\; \sum_{n=-\infty}^{\infty} u_{n,2}^{s} e^{i \kappa_n x_1 - i\zeta_n  x_2 } \quad \mbox{in} \; \Omega_2, \label{eq:rad_cond2}
\end{eqnarray}
in which the coefficients $u_{n,i}^s$ are complex amplitudes. In the above, the constants $\kappa_n$ and $\zeta_n$ are defined by
\begin{equation}
  \kappa_n=\kappa+\textstyle nb  \quad \mbox{and} \quad
\zeta_n = \zeta_n(k,\kappa)= \sqrt{k^2-\kappa_n^2}\,,
\end{equation}
where the domain of the analytic square root function is taken to be $\mathbb{C}\backslash\{-it: t\geq 0\}$, with $\sqrt{1}=1$.  With this choice of square root,
\begin{equation}
 \zeta_n(\kappa,k)   = \left\{
\begin{array}{lll}
\vspace*{5pt}
\sqrt{k^2-\kappa_n^2}\,,  & \mbox{if}\;\abs{\kappa_n} \leq k, \\
i\sqrt{\kappa_n^2-k^2}\,,  & \mbox{if}\;\abs{\kappa_n} \geq k. \\
\end{array}
\right. 
\end{equation}
The Rayleigh modes with $\abs{\kappa_n} < k$ are propagating, and the modes with $\abs{\kappa_n} > k$ are evanescent.  
The case of $\zeta_n=0$ corresonds to the so-called Rayleigh cut-off frequency, which is a transition frequency when
certain propagating modes convert to evanescent modes or vice versa \cite{rayleigh}.
Due to the redistribution of energy among different diffraction orders,
the tranmission spectral line is usually nonsmooth near such frequency as depicted in
Figure \ref{fig:trans_fano1} - \ref{fig:trans_fano3}, 
We won't be concerned with it in this paper.  Correspondingly, for given $\kappa$ and $k$, we define two sets
\begin{equation}\label{eq:Z1Z2}
\mathbb{Z}_1(\kappa,k) := \{ n \in \mathbb{Z} \,;\,  |\kappa_n|  < k  \} 
  \quad \mbox{and} \quad \mathbb{Z}_2(\kappa,k)  := \{ n \in \mathbb{Z} \,;\,   |\kappa_n|  > k    \} . 
\end{equation}

Due to the quasi-periodicity of the scattering problem, let us restrict the Bloch wave number $\kappa$ to the first Brillouin zone $\kappa\in(-b/2,b/2]$.
On the $\kappa$-$k$ plane, we introduce the diamond-shaped region
\begin{equation}\label{eq:D1}
  D_1 = \left\{ (\kappa,k) ; \, |\kappa|<b/2,\,  |\kappa|<k< b-\kappa \right\}.
\end{equation}
This is the parameter regime in which $\mathbb{Z}_1(\kappa,k)=\{0\}$ and there is exactly one diffraction mode.
As such we call $D_1$ the first radiation (diffraction) continuum.
Let 
\begin{equation}\label{eq:D2}
  D_2 = \left\{ (\kappa,k) ; \, |\kappa|<b/2,\,  k > b-\kappa \right\}
\end{equation}
be the region above $D_1$ where more than one diffraction mode are present. It is known that $\bar{D_1} \cup \bar{D_2}$ corresponds to the continuous spectrum
of the quasi-periodic scattering operator for all $\kappa$ in the first Brillouin zone \cite{bonnet_starling94}.

\subsection{BICs and Fano resonance with weakly coupled resonators}\label{sec:Fano_structure1}
Recall that the distance between the two slit holes in each period is $(\ell-1)\varepsilon$.  When they are not close to each other and there holds $\ell\varepsilon=O(1)$, the interaction between the two resonators are weak,
and we call them weakly coupled. Mathematically, the wave field generated from the interaction between the two resonators appear in the high-order terms of the total diffracted field.
This is in contrast to the case when $\ell=O(1)$ and the distance between the two slits holes is on the order $\varepsilon$.
In such a configuration, the wave field generated from the interaction between the two resonators is dominant in the total diffracted field. This would become clear
when we express the solution of the scattering problem as layer potentials in the next section, see, for instance, the expansions \eqref{eq:decomp_op}.

For weakly coupled resonators, due to the symmetry of the grating structure, embedded eigenvalues and BICs exist in the first diffraction continuum $D_1$,
and Fano resonance appears when the Bloch wavenumber $\kappa$ is perturbed from $0$ to break the symmetry. 
This was recently investigated in \cite{lin_shipman_zhang}, where
it is shown that two group of resonances admit the following asymptotic expansion if $|\kappa|\ll 1$ and $(\kappa,m\pi)\in D_1$:
\begin{equation}
k_m^{(j)}(\kappa) =  m \pi + 2m \varepsilon \ln\varepsilon + c_m^{(j)}(\kappa) \varepsilon + O(\varepsilon^2\ln^2\varepsilon), \quad j=1,2,
\end{equation}
where for fixed $\kappa$, $c_m^{(1)}$ and $c_m^{(2)}$ are constants independent of $\varepsilon$. 
Furthermore, there holds $\Im \, k_m^{(1)} = O(\varepsilon)$ and  $\Im \, k_m^{(2)} = O(\kappa^2 \varepsilon)$.

\begin{figure}[!htbp]
\begin{center}
\includegraphics[height=4.5cm, width=5cm]{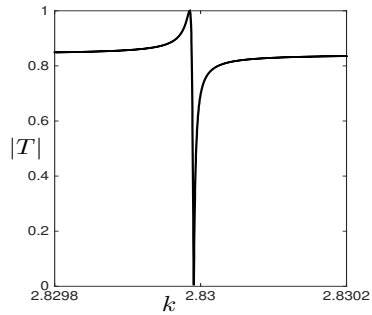}
\caption{Transmission $|T|$ when $d=1$, $\varepsilon=0.05$, $\kappa=0.1$, and the distance between two slits is $0.4$ within one period.
Fano resonance occurs near $k=2.83$. 
}\label{fig:trans_weak_coup}
\end{center}
\end{figure}

When $\kappa=0$, the imaginary parts of $k_m^{(2)}$ vanish, and the second group of resonances $k_m^{(2)}$ are
real-valued eigenvalues that are embedded  in the continuous spectrum $[0, +\infty)$.
If the Bloch wavenumber is perturbed with $ 0<|\kappa| \ll 1$,  then $k_m^{(2)}$ become complex-valued resonances with much smaller imaginary parts than $k_m^{(1)}$.
Therefore, the interaction of the two resonant modes with different strength yields the Fano transmission anomaly as shown in Figure \ref{fig:trans_weak_coup}. 
The readers are referred to \cite{lin_shipman_zhang} for a rigorous proof the appearance of the Fano type transmission signal,
and the studies for the associated field amplification at the resonant frequencies.
We emphasize that the above Fano resonance is induced by BICs, and it
would not occur in the higher diffraction continuum $D_2$, 
for which $k_m^{(1)}$ and $k_m^{(2)}$ attain the same order $O(\varepsilon)$ for the imaginary part, and consequently, the same order of resonant strength.

\subsection{Fano resonance with strongly coupled resonators}
When $\ell=O(1)$ and the two slit holes in each period are strongly coupled, there will still be two groups of resonances $k_m^{(1)}$ and $k_m^{(2)}$ close to each other,
and their asymptotic analysis will be derived in Theorem \ref{thm:asym_eig_perturb}.
For the imaginary parts of the resonances, there holds $\Im \, k_m^{(1)} = O(\varepsilon)$ if $(\kappa, m\pi) \in D_1 \cup D_2$.
On the other hand, $\Im \, k_m^{(2)} =  O(\kappa\varepsilon^2)$ if $(\kappa, m\pi) \in D_1$ and
$\Im \, k_m^{(2)}=O(\varepsilon^2)$ if $(\kappa, m\pi) \in D_2$ respectively (cf. Proposition \ref{prop:imag_res}).
Namely,  when $\kappa=0$, $k_m^{(2)}$ are still embedded eigenvalues in the first diffraction continuum, and similar to the weakly coupled case,
the perturbation of the Bloch wavenumber would give rise to Fano resonance.
However, above the diffraction continuum, the imaginary part of $k_m^{(2)}$ is consistently on the oder $\varepsilon^2$, versus order $\varepsilon$
for the imaginary part of $k_m^{(1)}$. 
As discussed in Section \ref{sec:anomaly},
the interference of these two resonant scattering, with one broad and the other narrow, will generate the Fano type transmission anomaly in a stable manner in this region.
Therefore,  compared to the Fano resonance induced by BICs,
the approach of strongly coupled resonators allows for the appearance of Fano resonance in a much broader frequency band.

The rest of the paper is organized as follows. We present an equivalent integral equation formulation for the scattering problem \eqref{eq:Helmholtz}--\eqref{eq:rad_cond2} and carry out the asymptotic expansions of the integral operators in Section~\ref{sec:bie}.
The asymptotic analysis of two groups of resonances is established in Section~\ref{sec:res}.
Based on the analysis of the diffracted field, the transmission anomaly is investigated in details in 
Section~\ref{sec:anomaly}.

\section{Boundary-integral formulation and asymptotic analysis}\label{sec:bie}
For clarity of exposition, we shall assume that $\ell = 2$ in the rest of the paper, and all the calculations can be carried over for $\ell \neq 2$ only with slight modification.
The scattering problem \eqref{eq:Helmholtz}--\eqref{eq:rad_cond2} can be formulated as a system of boundary-integral equations.
We refer to \cite{habib-book} for a systematic applications of the boundary-integral-equation method to various scattering problems in photonics and phononics.
This part is standard and is derived in Section 2 of \cite{lin_shipman_zhang}.  Here we only collect several key formulations.
For each $\kappa\in(-b/2,b/2]$, 
let $g_1(x,y)$ and $g_2(x,y)$ be the $d$-periodic Green's function in the semi-infinite domain $\Omega_1$ and $\Omega_2$
with the zero Neumann boundary condition along its boundary, respectively.
It is known that $g_1(x,y)$ and $g_2(x,y)$  are given by the sum of 
the free-space periodic Green's function $g(x,y)$ and $g(x',y)$, where $x'$ is the reflection point of $x$ with respect to the boundary of $\Omega_1$ and $\Omega_2$ respectively.
Let $g_\varepsilon^{\pm}(x,y)$ be the Green's function in the slit $S_{\varepsilon}^{0,\pm}$ with the zero Neumann boundary condition along the boundaries. We may express 
$g_\varepsilon^{\pm}(x,y) = g_\varepsilon^{0}(x_1 \mp \varepsilon,x_2; y_1 \mp \varepsilon, y_2)  $,
where $g_\varepsilon^{0}(x,y)$ is the Green's function in the rectangular domain $(-\frac{\varepsilon}{2},\frac{\varepsilon}{2})\times(0,1)$.

Let $\Omega^{(0)}:=\{ x \in \mathbf{R}^2 \; | \;   -\frac{d}{2}<x_1< \frac{d}{2} \}$ be reference period.
By applying the Green's theorem, one obtains
\begin{eqnarray*}
u_\varepsilon(x) &=& \int_{\Gamma^+_{1,\varepsilon} \cup \Gamma^-_{1,\varepsilon}} g_1(x,y) \dfrac{\partial u_\varepsilon(y)}{\partial\nu} ds_y + u^\mathrm{inc}(x)+ u^\mathrm{refl}(x)  \quad \mbox{for} \;\; x\in\Omega^{(0)} \cap \Omega_1, \\
u_\varepsilon(x) &=& \int_{\Gamma^+_{2,\varepsilon} \cup \Gamma^-_{2,\varepsilon}} g_2(x,y) \dfrac{\partial u_\varepsilon(y)}{\partial\nu} ds_y  \quad \mbox{for} \;\; x\in\Omega^{(0)} \cap \Omega_2, \\
u_\varepsilon(x) &=& -\int_{\Gamma^-_{1,\varepsilon} \cup \Gamma^-_{2,\varepsilon}} g_\varepsilon^{\mathrm{i},\pm}(x,y) \dfrac{\partial u_\varepsilon(y)}{\partial\nu} ds_y \quad \mbox{for} \;\; x\in S_\varepsilon^{0,\pm}.
\end{eqnarray*}
In the above, $\nu$ denotes the unit outward normal direction pointing toward the domain $\Omega_1$ or $\Omega_2$,
and $u^\mathrm{refl}(x)=e^{i (\kappa x_1 + \zeta_0 (x_2-1))}$ is the reflected field of the ground plane $\{x_2=1\}$ without the slits.

Taking the limit of the layer potentials to the slit apertures $\Gamma^\pm_{1,\varepsilon}$ and $\Gamma^\pm_{2,\varepsilon}$ from the exterior and interior of the slits, 
and imposing the continuity condition for $u_\varepsilon$ along them leads to the following system of four integral equations:
\begin{equation}\label{eq:scattering2}
\left\{
\begin{array}{l}
\displaystyle{\int_{\Gamma^+_{1,\varepsilon} \cup \Gamma^-_{1,\varepsilon}}  g_1(x,y)\dfrac{\partial u_\varepsilon(y)}{\partial \nu} ds_y 
+\int_{\Gamma^-_{1,\varepsilon} \cup \Gamma^-_{2,\varepsilon}} g_\varepsilon^{\mp}(x,y) \dfrac{\partial u_\varepsilon(y)}{\partial \nu} ds_y }
= - (u^\mathrm{inc}(x)+u^\mathrm{refl}(x)), \quad \mbox{for} \,\, x\in\Gamma^\mp_{1,\varepsilon}, \\ \\ 
\displaystyle{\int_{\Gamma^+_{2,\varepsilon} \cup \Gamma^-_{2,\varepsilon}}  g_2(x,y)\dfrac{\partial u_\varepsilon(y)}{\partial \nu} ds_y
+\int_{\Gamma^-_{1,\varepsilon} \cup \Gamma^-_{2,\varepsilon} } g_\varepsilon^{\mp}(x,y) \dfrac{\partial u_\varepsilon(y)}{\partial \nu} ds_y
=0},  \quad \mbox{for} \,\, x\in\Gamma^\mp_{2,\varepsilon}.
\end{array}
\right.
\end{equation}
Note that on the slit apertures, there holds $x=x'$, and the periodic Green's functions satisfy
\begin{equation}\label{eq:g12}
g_1(x,y)= g_2(x,y) = 2 g(x_1,0;y_1,0).
\end{equation}
On the other hand, over the slit apertures, the slit Green's functions can be expressed as
\begin{eqnarray}
 g_\varepsilon^{\pm}(x_1,1;y_1,1) = g_\varepsilon^{\pm}(x_1,0;y_1,0) 
&=&  g_\varepsilon^{0}(x_1 \mp \varepsilon, 0; y_1 \mp \varepsilon, 0),  \label{eq:g_slit1}  \\
 g_\varepsilon^{\pm}(x_1,1;y_1,0) = g_\varepsilon^{\pm}(x_1,0;y_1,1) 
&=&  g_\varepsilon^{0}(x_1 \mp \varepsilon, 1; y_1 \mp \varepsilon, 0).  \label{eq:g_slit2}
\end{eqnarray}

We rescale the slit aperture to the unit interval $I : = (-\frac{1}{2},\frac{1}{2})$ by letting
$$ x_1 =  \varepsilon  (X \pm 1)  \quad \mbox{for} \; (x_1,1) \in \Gamma_{1,\varepsilon}^\pm \;\; \mbox{and} \; (x_1,0) \in \Gamma_{2,\varepsilon}^\pm.$$
From the Rayleigh-Bloch expansion of the periodic Green's function (cf.~\cite{linton98}),
the rescaled periodic Green's function over the slit apertures $\Gamma^\pm_{1,\varepsilon}$ and $\Gamma^\pm_{2,\varepsilon}$ are given explictly by
$$ G_\varepsilon^\mathrm{e}(X, Y) := 2 g(\varepsilon X, 0; \varepsilon Y, 0) = -\dfrac{i}{d} \sum_{n=-\infty}^{\infty} \dfrac{1}{\zeta_n(\kappa,k)   } e^{i \kappa_n\varepsilon(X-Y)}.   $$
Similarly, the slit Green's function is rescaled and expressed explicitly as
 \begin{eqnarray*}
 G_\varepsilon^\mathrm{i}(X, Y) &:=&  g_\varepsilon^{0}( \varepsilon X, 0; \varepsilon Y, 0 ) 
 = \frac{1}{\varepsilon} \sum_{m,n=0}^\infty c_{mn}a_{mn}\cos\left(m\pi (X+\textstyle\frac{1}{2})\right ) \cos\left(m\pi (Y+\textstyle\frac{1}{2})\right ); \\ 
\tilde  G_\varepsilon^\mathrm{i}(X, Y) &:=&  g_\varepsilon^{0}(\varepsilon X, 1; \varepsilon Y, 0) 
 = \frac{1}{\varepsilon} \sum_{m,n=0}^\infty(-1)^n c_{mn}a_{mn} \cos\left(m\pi (X+\textstyle\frac{1}{2})\right ) \cos\left(m\pi (Y+\textstyle\frac{1}{2})\right ),
\end{eqnarray*}
in which the coefficients 
\begin{equation*}
c_{mn}=[k^2-(m\pi/\varepsilon)^2 - (n\pi)^2]^{-1}
\quad \mbox{and} \quad
a_{mn} = \left\{
\begin{array}{llll}
1  & m=n=0, \\
2  & m=0, n\ge 1 \quad \mbox{or} \quad n=0, m\ge 1, \\
4  & m\ge 1, n \ge 1.
\end{array}
\right.
\end{equation*}
Correspondingly, let us introduce the rescaled boundary-integral operators for $X\in I$:
\begin{eqnarray}\label{eq:op_T}
 && [T^\mathrm{e} \varphi](X) = \int_I G_\varepsilon^\mathrm{e}(X, Y)  \varphi(Y) dY,  \quad
 [T^{\mathrm{e},\pm} \varphi](X) = \int_I G_\varepsilon^\mathrm{e}(X, Y \mp 2) \varphi(Y) dY   \label{eq:op_Te} \\
 && [T^\mathrm{i}  \varphi] (X) = \int_I  G_\varepsilon^\mathrm{i}(X, Y) \varphi(Y) dY, \quad 
    [\tilde T^\mathrm{i}  \varphi] (X) = \int_I  \tilde G_\varepsilon^\mathrm{i}(X, Y) \varphi(Y) dY. \label{eq:op_Ti} 
\end{eqnarray}

If one defines the following quantities in the scaled interval,
\begin{eqnarray*}
&& \varphi_1^\pm(X):= \dfrac{\partial u_\varepsilon}{\partial\nu}( \varepsilon (X\pm 1), 1),  \varphi_2^\pm(X):= \dfrac{\partial u_\varepsilon}{\partial\nu}( \varepsilon (X\pm 1), 0),   \\
&& f^\pm(X):= -\frac{1}{2} (u^\mathrm{inc}+u^\mathrm{refl})(\varepsilon (X \pm 1), 1) = -e^{i \kappa \varepsilon (X\pm 1 )},
\end{eqnarray*}
then in view of \eqref{eq:g12} - \eqref{eq:g_slit2}, the system (\ref{eq:scattering2}) can be recast as 
$\mathbb{T}\boldsymbol{\varphi}=\varepsilon^{-1}\mathbf{f}$, in which
\begin{equation}\label{eq:system_BIE}
\mathbb{T}=\left[
\begin{array}{cccc}
T^\mathrm{e}+T^\mathrm{i}    &  T^{\mathrm{e},-}   & \tilde T^\mathrm{i} & 0  \\
 T^{\mathrm{e},+}     &  T^\mathrm{e}+T^\mathrm{i} & 0 & \tilde T^\mathrm{i}     \\
  \tilde T^\mathrm{i}  & 0  & T^\mathrm{e}+T^\mathrm{i}    &  T^{\mathrm{e},-}  \\
   0  & \tilde T^\mathrm{i}  & T^{\mathrm{e},+}  & T^\mathrm{e}+T^\mathrm{i}   
\end{array}
\right], \quad
\boldsymbol{\varphi} =
\left[
\begin{array}{cccc}
\varphi_1^-    \\
\varphi_1^+   \\
\varphi_2^-   \\
\varphi_2^+
\end{array}
\right], \quad
\mathbf{f}=\left[
\begin{array}{cccc}
2f^-  \\
2f^+  \\
0   \\
0
\end{array}
\right].
\end{equation}

\bigskip

For each fixed $\kappa\in(-b/2,b/2]$, we derive the asymptotic expansion of the boundary integral operators $T^e$, $T^{e,\pm}$,
$T^i$, and $\tilde{T}^i$ when the slit size $\varepsilon \ll 1$. Here and henceforth, we assume that $k=O(1)$ and $k$ is away from the Rayleigh cut-off frequencies $\kappa+nb$ for all integers $n$. By doing this we exclude the scenario where the Green function $G_\varepsilon^e(X,Y)$ is not defined because of the vanishing value for certain $\zeta_n$.

It can be shown that the kernel $G_\varepsilon^\mathrm{e}(X,Y)$ admits the following asymptotic expansion: 
\begin{equation}\label{Ge_exp}
G_\varepsilon^\mathrm{e}(X, Y)=\beta_\mathrm{e}(\kappa,k,\varepsilon) + \dfrac{1}{\pi}  \ln |X-Y| +  r_\mathrm{e}(\kappa,\varepsilon; X,Y), 
\end{equation}
where
\begin{equation}\label{beta_e}
\beta_\mathrm{e}(\kappa,k, \varepsilon)= \dfrac{1}{\pi} \ln (\varepsilon b)  + \sum_{n\neq 0} \left( \dfrac{1}{2\pi} \dfrac{1}{|n|} - \dfrac{i}{d}  \dfrac{1}{\zeta_n(\kappa,k)}\right)
  - \dfrac{i}{d}  \dfrac{1}{\zeta_0(\kappa,k)}.
\end{equation}
The high-order term satisfies
\begin{equation}\label{eq:re_D1}
r_\mathrm{e}(\kappa,\varepsilon; X,Y)=r_\mathrm{e}(0,\varepsilon; X,Y)+O(\kappa\varepsilon) \quad \mbox{if} \; (\kappa,k) \in D_1,
\end{equation}
where $r_\mathrm{e}(0,\varepsilon; X,Y)=O(k^2\varepsilon^2|\ln\varepsilon|)$ is real-valued. On the other hand,
\begin{equation}\label{eq:re_D2}
r_\mathrm{e}(\kappa,\varepsilon; X,Y)=O(k\varepsilon) \quad \mbox{if} \; (\kappa,k) \in D_2.
\end{equation}
The regions $D_1$ and $D_2$ are defined by \eqref{eq:D1} and \eqref{eq:D2}.
The above expansion for $G_\varepsilon^\mathrm{e}(X, Y)$ is derived in the appendix.
The asymptotic expansions of $G_\varepsilon^\mathrm{i}(X, Y)$ and $\tilde G_\varepsilon^\mathrm{i}(X, Y)$ are given by (cf. Lemma 3.1 of \cite{lin_zhang17})
\begin{eqnarray}\label{eq:Gi_asy}
G_\varepsilon^\mathrm{i}(X, Y) 
 &=& \beta_\mathrm{i}(k ,\varepsilon) + \dfrac{1}{\pi} \left[ \ln \abs{\sin \left(\frac{\pi(X-Y)}{2}\right)} + \ln \abs{\sin \left(\frac{\pi(X+Y+1)}{2}\right)} \right] \nonumber \\
 && +\; r_{\mathrm{i}}(\varepsilon; X,Y);     \label{Gi_exp}  \\
\tilde G_\varepsilon^\mathrm{i}(X, Y) &= &  \tilde \beta(k,\varepsilon) +  \tilde r_{\mathrm{i}}(\varepsilon; X,Y).
 \label{tGi_exp} 
\end{eqnarray}
In the above, 
\begin{equation}\label{beta_i}
 \beta_\mathrm{i}(k, \varepsilon)= \dfrac{1}{\varepsilon\, k\tan k } +  \dfrac{2\ln 2}{\pi},  \quad\quad \tilde \beta(k,\varepsilon) = \dfrac{1}{ \varepsilon\,k\sin k},  
\end{equation}
and the real-valued high-order terms $r_{\mathrm{i}} = O(\varepsilon^2)$ and $\tilde r_{\mathrm{i}} = O(e^{-1/\varepsilon})$.

Let $P$ be the projection operator on the interval $I$ defined by $ P \varphi(X) = \langle \varphi, 1 \rangle 1$, where $1$ is a function defined on the interval $I$ and is equal to one therein. 
Define the functions spaces 
 $$V_1 = \tilde H^{-\frac{1}{2}}(I):=\{ u=U|_{I}  \;\big{|}\;  U \in H^{-1/2}(\mathbf{R})\,\,  \mbox{and} \,\, supp \,U \subset \bar I   \} \quad \mbox{and} \quad V_2 =  H^{\frac{1}{2}}(I). $$
We introduce the following integral operators from $V_1$ to $V_2$:
\begin{eqnarray*}\label{eq:op_Ks}
 && [S \varphi](X) = \dfrac{1}{\pi}  \int_{I} \ln \abs{ (X-Y) \sin \left(\frac{\pi(X-Y)}{2}\right) \sin \left(\frac{\pi(X+Y+1)}{2}\right)} \, \varphi(Y) \, dY; \label{eq:op_K} \\
 && [S^{\pm} \varphi](X) = \dfrac{1}{\pi} \int_I  \ln |X-Y \pm 2|  \, \varphi(Y) \, dY; \label{eq:op_K_12} \\
 && [S^{\infty} \varphi](X) = \int_I  \Big[  r_\mathrm{e}(X,Y) + r_{\mathrm{i}}(X,Y) \Big]\,  \varphi(Y) \, dY; \label{eq:op_K_infty} \\
 && [S^{\infty,\pm} \varphi](X) =  \int_I   r_e(X,Y\mp 2) \, \varphi(Y) \, dY; \label{eq:op_K_12} \\
 && [\tilde S^{\infty} \varphi](X) = \int_I   \tilde r_\mathrm{i} (X,Y) \, \varphi(Y) \, dY.  \label{eq:op_K_21}
\end{eqnarray*}
For brevity, in the above we neglect the explicit dependence of the kernels $r_\mathrm{e}$, $r_\mathrm{i}$ and $\tilde r_\mathrm{i}$ on $\kappa$ and $\varepsilon$.
Then from the asymptotic expansion of the Green functions,
the integral operators $T^\mathrm{e}+T^\mathrm{i} $, $T^{\mathrm{e},\pm}$, and $\tilde T^\mathrm{i}$ can be decomposed as
\begin{equation} \label{eq:decomp_op}
T^\mathrm{e}+T^\mathrm{i}  = (\beta_\mathrm{e} + \beta_\mathrm{i})  P + S + S^\infty, \quad T^{\mathrm{e},\pm} =  \beta_\mathrm{e} P + S^{\pm}+S^{\infty,\pm}, 
\quad \tilde T^\mathrm{i} = \tilde\beta P + \tilde S^{\infty},
\end{equation}
where $S_\infty$, $S^{\infty,\pm}_\kappa$ and $S^{\infty}$ are bounded from $V_1$ to $V_2$, with operator norm
\begin{equation}\label{eq:S_infty_norm1}
\| S^{\infty}\|  \lesssim \varepsilon, \quad  \| S^{\infty,\pm} \|  \lesssim \varepsilon,
\quad \mbox{and} \quad \| \tilde S^{\infty}\|  \lesssim e^{-1/\varepsilon}.
\end{equation}

If $(\kappa,k) \in D_1$, in view of the relation \eqref{eq:re_D1} and the estimation $r_\mathrm{e}(0,\varepsilon; X,Y)=O(k^2\varepsilon^2|\ln\varepsilon|)$,
it follows that $r_\mathrm{e}(\kappa,\varepsilon; X,Y) = O(k^2\varepsilon^2 |\ln \varepsilon|)  + O(|\kappa|\varepsilon)$.
In addition, note that $r_{\mathrm{i}}(X,Y) = O(\varepsilon^2)$. Thus we have the following lemma for $S^{\infty}$ and $S^{\infty,\pm}$. 

\begin{lem}\label{lem:hot}
If $(\kappa,k) \in D_1$, then there holds that
\begin{equation}\label{eq:S_infty_norm2}
\| S^{\infty}\|  \lesssim k^2\varepsilon^2 |\ln \varepsilon|  + |\kappa|\varepsilon \quad \mbox{and} \quad  \| S^{\infty,\pm} \|  \lesssim k^2\varepsilon^2 |\ln \varepsilon|  + |\kappa|\varepsilon.
\end{equation} 
\end{lem}

\section{The resonances for the scattering problem} \label{sec:res}
\subsection{The resonance condition}\label{sec:res_cond}
For each $\kappa \in (-b/2,b/2]$, we solve for complex-valued resonances $k$, whose real part lies in
the continuous spectrum of the quasi-periodic scattering operator.
The corresponding quasi-mode grows exponentially away from the grating.
From the discussions in Section \ref{sec:bie}, we see that the resonances are the complex characteristic values $k$ of the system of integral equations
\eqref{eq:scattering2}, or equivalently, the complex characteristic frequencies of the operator $\mathbb{T}$ in the scaled intervals for which the system $\mathbb{T}\, \boldsymbol{\varphi}=0$ admits a nonzero solution.

Let $\sigma(\mathbb{T})$ be the set of complex characteristic frequencies $k$ of $\mathbb{T}$.  Note that
\begin{equation*}
\mathbb{T}=\left[
\begin{array}{cc}
\hat{\mathbb{T}}  & \tilde{\mathbb{T}}  \\
\tilde{\mathbb{T}}  & \hat{\mathbb{T}}
\end{array}
\right], \quad\mbox{where} \;
\hat{\mathbb{T}}=\left[
\begin{array}{cc}
T^\mathrm{e}+T^\mathrm{i}    &  T^{\mathrm{e},-}   \\
 T^{\mathrm{e},+}     &  T^\mathrm{e}+T^\mathrm{i}  
\end{array}
\right]
\; \mbox{and} \;
\tilde{\mathbb{T}}=\left[
\begin{array}{cc}
 \tilde T^\mathrm{i} & 0  \\
 0 & \tilde T^\mathrm{i}  \\ 
\end{array}
\right].
\end{equation*}
One can decompose the function space $(V_1)^4$ as $(V_1)^4=V_\mathrm{even} \oplus V_\mathrm{odd}$,
where $V_\mathrm{even}$ and $V_\mathrm{odd}$ are invariant subspaces for~$\mathbb{T}$ given by
$ V_\mathrm{even} = \{ \,[ \varphi_-, \varphi_+,  \varphi_-, \varphi_+ ]^T \,;\, \varphi_\pm \in V_1 \}  $
and
$ V_\mathrm{odd} = \{ \,[ \varphi_-, \varphi_+,  -\varphi_-, -\varphi_+ ]^T \,;\, \varphi_\pm \in V_1 \} $.
As such, it follows that
$$\sigma(\mathbb{T})= \sigma(\mathbb{T}|_{V_\mathrm{even}})\, \cup \, \sigma(\mathbb{T}|_{V_\mathrm{odd}})=
\sigma\big(\mathbb{T}_+\big)\, \cup \, \sigma\big(\mathbb{T}_-\big),$$
where 
\begin{equation}\label{eq:T_pm}
\mathbb{T}_+=\hat{\mathbb{T}}+\tilde{\mathbb{T}} \quad \mbox{and} \quad \mathbb{T}_-=\hat{\mathbb{T}}-\tilde{\mathbb{T}}.
\end{equation}
Hence, we shall derive the characteristic values of the operators $\mathbb{T}_+$ and $\mathbb{T}_-$ to obtain $\sigma(\mathbb{T})$.
This can be reduced to solving for the roots of certain nonlinear functions as shown it what follows.

Let
\begin{equation}\label{eq:beta}
\beta(\kappa,k, \varepsilon)=\beta_\mathrm{e}(\kappa,k, \varepsilon) + \beta_\mathrm{i}(k, \varepsilon) - \beta_0,
\end{equation}
where the quantity $\beta_\mathrm{e}$ and $\beta_\mathrm{i}$ is given in \eqref{beta_e} and \eqref{beta_i} respectively,
$\beta_0$ is certain real constant independent of $\varepsilon$ and $k$ to be specified in Lemma \ref{lem:Q_hat}.
Introduce the operators
\begin{equation}\label{eq:op_PS}
\mathbb{P}_\pm= \left[
\begin{array}{cc}
(\beta\pm\tilde\beta) P    & \beta_\mathrm{e} P \\
\beta_\mathrm{e} P  & (\beta\pm\tilde\beta) P 
\end{array}
\right],
\;
\mathbb{S}=
\left[
\begin{array}{cc}
S+\beta_0 P   & S^{-}  \\
S^{+}  & S +\beta_0 P
\end{array}
\right],  \;
\mathbb{S}_\pm^\infty=
\left[
\begin{array}{cc}
S^\infty \pm \tilde S^{\infty}    & S^{\infty,-} \\
S^{\infty,+}  & S^\infty  \pm \tilde S^{\infty}
\end{array}
\right],
\end{equation}
In view of the decomposition for the integral operators in \eqref{eq:decomp_op}, it follows that
$$ \mathbb{T}_\pm = \mathbb{P}_\pm + \mathbb{S}  + \mathbb{S}_\pm^\infty =:  \mathbb{P} + \mathbb{L}_\pm, $$
and the homogeneous problem $\mathbb{T}_\pm\, \boldsymbol{\varphi}=0$ become
\begin{equation}\label{eq:scattering4}
(\mathbb{P}_\pm + \mathbb{L}_\pm) \boldsymbol{\varphi} =   \mathbf{0}, 
\end{equation}
where $\boldsymbol{\varphi} = [\, \varphi^-,  \varphi^+ \,] ^T$.

It is shown in \cite{eric10-2} that the operator $\mathbb{S}$ attains a bounded inverse from $(V_2)^2$ to $(V_1)^2$.
Consequently, $\mathbb{L}_\pm$ is invertible for sufficiently small $\varepsilon$ since $\| \mathbb{S}_\pm^{\infty}\|  \lesssim \varepsilon$. 
By applying $\mathbb{L}_\pm^{-1}$ on both sides of \eqref{eq:scattering4}, it can be calculated explicitly that
\begin{equation}\label{Op_eqns_1}
\left[ (\beta\pm\tilde \beta)  \langle  \boldsymbol{\varphi}, \mathbf{e}_1 \rangle  + 
 \beta_\mathrm{e} \langle  \boldsymbol{\varphi}, \mathbf{e}_2 \rangle \right] \, \mathbb{L}_\pm^{-1} \mathbf{e}_1  
 + \left[ \beta_\mathrm{e} \langle  \boldsymbol{\varphi}, \mathbf{e}_1 \rangle  + (\beta\pm\tilde \beta) \langle \boldsymbol{\varphi}, \mathbf{e}_2 \rangle \right] \,
  \mathbb{L}_\pm^{-1}  \mathbf{e}_2 + \boldsymbol{\varphi} =   \mathbf{0},
 \end{equation}
where $\{\mathbf{e}_j\}_{j=1}^2\in V_1 \times V_1 $  are given by $\mathbf{e}_1 = [1, 0]^T $ and $\mathbf{e}_2 = [0, 1]^T $. 
Taking the $L^2$-inner product of \eqref{Op_eqns_1} with $\mathbf{e}_1$ and $\mathbf{e}_2$ yields a $2 \times 2$ linear system
\begin{equation}\label{eq:linear_sys1}
\tilde{\mathbb{M}}_{\pm} \left[
 \begin{array}{llll}
\langle  \boldsymbol{\varphi}, \mathbf{e}_1 \rangle  \\
\langle  \boldsymbol{\varphi}, \mathbf{e}_2 \rangle
\end{array}
\right] =
\left[
\begin{array}{llll}
0   \\
0
\end{array}
\right],
\quad\mbox{where} \; \tilde{\mathbb{M}}_{\pm} (k; \kappa, \varepsilon):= \mathbb{Q}_{\pm}\mathbb{B}_\pm + \mathbb{I}.
\end{equation}
In the above, $\mathbb{I}$ is the identity matrix, and
\begin{equation}  \label{eq:matrix_QB}
\mathbb{Q}_{\pm}(k; \kappa, \varepsilon):=
\left[
\begin{array}{cc}
 \langle \mathbb{L}_{\pm}^{-1}  \mathbf{e}_1,  \mathbf{e}_1 \rangle  &  \langle \mathbb{L}_{\pm}^{-1}  \mathbf{e}_2,  \mathbf{e}_1\rangle \\
 \langle \mathbb{L}_{\pm}^{-1}  \mathbf{e}_1,  \mathbf{e}_2 \rangle  &  \langle \mathbb{L}_{\pm}^{-1}  \mathbf{e}_2,  \mathbf{e}_2 \rangle 
 \end{array}
\right], 
\quad 
\mathbb{B}_\pm (k; \kappa, \varepsilon) :=
 \left[
\begin{array}{cc}
 \beta \pm \tilde\beta &  \beta_\mathrm{e} \\
\beta_\mathrm{e} &  \beta \pm \tilde\beta
 \end{array}
\right].   
\end{equation}
Therefore, the characteristic values of the operator-valued function $\mathbb{T}_\pm(k;\kappa,\varepsilon)$ are the roots of 
$\tilde\lambda_{1,\pm}(k; \kappa, \varepsilon)$ and $\tilde\lambda_{2,\pm}(k; \kappa, \varepsilon)$, the eigenvalues of $\tilde{\mathbb{M}}_{\pm}$.
Since the leading-order of $\beta$ and $\beta_\mathrm{e}$ in $\varepsilon$ is $O(1/\varepsilon)$, we scale the matrix $\tilde{\mathbb{M}}_\pm$ by letting
$$
\mathbb{M}_{\pm} \,:=\, \varepsilon\, \tilde{\mathbb{M}}_{\pm}.$$
The eigenvalues of $\mathbb{M}_{\pm}$ are
\begin{equation}\label{eq:lambdas}
\lambda_{j,\pm}(k; \kappa, \varepsilon) \,:=\, \varepsilon \, \tilde\lambda_{j,\pm}(k; \kappa, \varepsilon), \quad j=1,2.
\end{equation}
The following proposition summarizes the resonance condition.

\begin{prop} \label{prop:res_cond}
For each $\kappa \in (-b/2,b/2]$, let
$\mathbb{C}_\kappa=\{ z \in \mathbb{C} \,;\,    \Re \, z > |\kappa| \} \subset \mathbb{C}$.
Then the complex-valued resonances of the scattering problem \eqref{eq:Helmholtz}--\eqref{eq:rad_cond2} 
in $\mathbb{C}_\kappa$ are the complex roots of the functions $\lambda_{j,\pm}(k; \kappa, \varepsilon)$, $(j=1,2)$, 
where $\lambda_{1,\pm}$ and $\lambda_{2,\pm}$ are eigenvalues of the matrix $\mathbb{M}_{\pm}$.
\end{prop}

\subsection{Asymptotic expansion of resonances}\label{sec:resonances}
From Proposition \ref{prop:res_cond}, one can solve for the roots of $\lambda_{j,\pm}(k; \kappa, \varepsilon)$ to obtain resonances. 
To this end, we will first approximate $\mathbb{M}_{\pm}$ by certain symmetric matrices whose eigenvalues are easier to analyze,
and then apply the perturbation theory to study the high-order terms of the roots. This method will be used several times throughout the rest of the paper. We introduce it in what follows.

\subsubsection{A perturbation theory for matrices}\label{sec:perturb_mat}
Given a matrix
\begin{equation*}
\mathbb{M}:= \varepsilon (\mathbb{Q}\mathbb{B} +  \mathbb{I}), 
\quad \mbox{where} \;
\mathbb{Q}:=\left[
\begin{array}{cc}
 q_{11}(\varepsilon)  & q_{12}(\varepsilon) \\
 q_{21}(\varepsilon) &  q_{22}(\varepsilon)
 \end{array}
\right], 
\quad
\mathbb{B}:=\left[
\begin{array}{cc}
 b(\varepsilon)  &  \tilde b(\varepsilon) \\
 \tilde b(\varepsilon) &  b(\varepsilon)
 \end{array}
\right].
\end{equation*}
Assume that there exits $\hat q_{ij}$ such that
$$ q_{ij}(\varepsilon)-\hat q_{ij}=O\big(\tau(\varepsilon)\big), \hat q_{11} = \hat q_{22} = q  \; \mbox{and}  \; \hat q_{12}=\hat q_{21} = \tilde q, $$
where $\tau(\varepsilon) \to 0$ as $\varepsilon\to 0$.
We introduce a symmetric matrix
\begin{equation}  \label{eq:M_hat}
\hat{\mathbb{M}}:=\varepsilon(\hat{\mathbb{Q}}\mathbb{B} + \mathbb{I})
\quad \mbox{where} \;
\hat{\mathbb{Q}} =
\left[
\begin{array}{cc}
 q  & \tilde q \\
 \tilde q &  q
 \end{array}
\right].
\end{equation}
 It can be calculated that the eigenvalues of $\hat{\mathbb{M}}$ are
\begin{equation}\label{eq:eigenval1_M_hat} 
\hat\lambda_{1} = \varepsilon+\varepsilon (b+\tilde b) (q + \tilde q), \quad
\hat\lambda_{2} = \varepsilon+\varepsilon (b- \tilde b) (q - \tilde q),
\end{equation}
and the corresponding eigenvectors are
\begin{equation}\label{eq:eigenvec_M_hat}
\hat v_{1} = [\; 1 \;\;\;\; 1 \; ] ^T, \quad
\hat v_{2} = [\; 1 \;\;\; -1 \; ] ^T.
\end{equation}
The sensitivity of eigenvalues of the matrix $\mathbb{M}$ with respect to 
the perturbation $\Delta  \mathbb{M} := \mathbb{M} - \hat{\mathbb{M}}$ is given in the following lemma.
\begin{lem}\label{lem:perturb_eig_M}
Let $\{\lambda_{j}\}_{j=1}^2$ and  $\{v_{j}\}_{j=1}^2$ be the eigenvalues and eigenvectors of $\mathbb{M}$.
Assume that $\hat{\mathbb{Q}}$ is invertible. If $\|\mathbb{Q}-\hat{\mathbb{Q}}\|=O(\tau)$ and $\|\mathbb{M}-\hat{\mathbb{M}}\| =O(\tau)$, then
\begin{equation}\label{eq:lambda_sensitivity}
\lambda_{j}   = (1+O(\tau) ) \cdot \hat\lambda_{j}  + O(\varepsilon \tau),   \quad
v_{j}  = \hat v_{j} +  O(\tau),  \quad j=1,2.
\end{equation}
\end{lem}

\noindent\textbf{Proof.}   Let  $\Delta\lambda_{j}=\lambda_{j} - \hat\lambda_{j}$  and $\Delta v_{j} =v_{j} - \hat v_{j}$ be the perturbation of the eigenvalues and eigenvectors.
The sensitivity of eigenvectors is obvious. To investigate the sensitivity of eigenvalues,
we expand the relation $\mathbb{M}  v_{j} = \lambda_{j} v_{j}$ to obtain
$$\hat \lambda_{j} \cdot \Delta v_{j} + \Delta\lambda_{j} \cdot \hat v_{j}  \;=\; \hat{\mathbb{M}} \, \Delta v_{j} +  \Delta  \mathbb{M} \, \hat v_{j} + O(\tau^2).  $$
Multiplying by $\hat v_{j}^T$ yields
\begin{equation}\label{eq:lambda_perturb1}
 2\Delta\lambda_{j}  \;=\; \hat v_{j}^T \, \Delta \mathbb{M}  \, \hat v_{j}  + O(\tau^2).
\end{equation}
Using the relation $\hat{\mathbb{M}}\hat v_{j} = \hat \lambda_{j} \hat v_{j} $,
one can express $\Delta \mathbb{M} \, \hat v_{j}$ as 
\begin{equation}\label{eq:lambda_perturb2}
\Delta \mathbb{M} \, \hat v_{j} =
\varepsilon(\mathbb{Q}-\hat{\mathbb{Q}}) \mathbb{B} \hat v_{j} =
 (\hat\lambda_{j} -\varepsilon) 
(\mathbb{Q}-\hat{\mathbb{Q}}) \hat{\mathbb{Q}}^{-1}
\hat v_{j},
\end{equation}
and the assertion follows by combining \eqref{eq:lambda_perturb1} and \eqref{eq:lambda_perturb2}.
\qed

\subsubsection{Asymptotic expansions of resonances}
We apply the perturbation method in the previous section to derive the asymptotic expansions for the roots of the functions $\lambda_{j,\pm}(k; \kappa, \varepsilon)$ $(j=1,2)$.
Recall that $\lambda_{1,\pm}$ and $\lambda_{2,\pm}$ are eigenvalues of the matrix $\mathbb{M}_{\pm}$.
We first introduce the following approximation for $\mathbb{M}_{\pm}$:
\begin{equation}  \label{eq:matrix-M_hat}
\hat{\mathbb{M}}_{\pm}:=\varepsilon(\hat{\mathbb{Q}}\mathbb{B}_\pm + \mathbb{I})
\quad \mbox{where} \;
\hat{\mathbb{Q}}
=
\left[
\begin{array}{cc}
 \langle \mathbb{S}^{-1}  \mathbf{e}_1,  \mathbf{e}_1 \rangle  &  \langle \mathbb{S}^{-1}  \mathbf{e}_2,  \mathbf{e}_1\rangle \\
 \langle \mathbb{S}^{-1}  \mathbf{e}_1,  \mathbf{e}_2 \rangle  &  \langle \mathbb{S}^{-1}  \mathbf{e}_2,  \mathbf{e}_2 \rangle 
 \end{array}
\right].
\end{equation}
In the above, the operator $\mathbb{S}$ is defined by \eqref{eq:op_PS} and the matrix $\mathbb{B}_\pm$ by \eqref{eq:matrix_QB}.

\begin{lem}\label{lem:Q_hat}
The matrix $\hat{\mathbb{Q}}$ is symmetric and there exists real constant $\beta_0$ such that $\hat{\mathbb{Q}}$ is invertible.
\end{lem}

\noindent\textbf{Proof.}
Let $\boldsymbol{\varphi}=(\varphi_1, \varphi_2)^T$ and $\tilde{\boldsymbol{\varphi}}=(\tilde{\varphi}_1, \tilde{\varphi}_2)^T$ satisfy
 $\mathbb{S}\boldsymbol{\varphi} =\mathbf{e}_1$ and $\mathbb{S}\tilde{\boldsymbol{\varphi}} =\mathbf{e}_2 $. 
By a direct calculation, it is seen that
 $$ (\hat S-S^{-} \hat S^{-1} S^{+})\varphi_1 = 1, \quad (\hat S-S^{+} \hat S^{-1} S^{-})\tilde{\varphi}_2 = 1,   $$
where $\hat S:= S+\beta_0 P$. From the definition of the operators $S$ and $S^{\pm}$, we see that for a function $\varphi$ defined over the interval $I$, there holds
$$ [S\tilde{\varphi}](X) = [S\varphi](-X),  \; [S^{-}\tilde{\varphi}](X) = [S^{+}\varphi](-X), \; [S^{+}\tilde{\varphi}](X) = [S^{-}\varphi](-X), $$
where $\tilde{\varphi}(X)=\varphi(-X)$.
Therefore, from the uniquess of the solution for the integral equation,  $\tilde{\varphi}_2(X)=\varphi_1(-X)$ holds and consequently 
$\langle \mathbb{S}^{-1}  \mathbf{e}_1,  \mathbf{e}_1 \rangle = \langle \mathbb{S}^{-1}  \mathbf{e}_2,  \mathbf{e}_2 \rangle$.
Similarly, it can be shown that  $\tilde{\varphi}_1(X)=\varphi_2(-X)$ and $\langle \mathbb{S}^{-1}  \mathbf{e}_1,  \mathbf{e}_2 \rangle = \langle \mathbb{S}^{-1}  \mathbf{e}_2,  \mathbf{e}_1 \rangle$. This proves the symmetry of the matrix $\hat{\mathbb{Q}}$.

We follow the lines in \cite{eric10-2} to show the invertibility of $\hat{\mathbb{Q}}$.
Introduce the operator $\mathbb{T}: (V_1)^2 \times \mathbb{R}^2 \to (V_2)^2 \times \mathbb{R}^2$ such that 
$\mathbb{T}(\mathbf{g}, \mathbf{u}):=\Big((\mathbb{S}-\beta_0 P)\mathbf{g}-\mathbf{u}, \displaystyle{\int_I \mathbf{g}} \, dx \Big)$,
 then it can be shown that $\mathbb{T}$ is invertible (Lemma 4.4 of \cite{eric10-2}). Let $(\boldsymbol{\psi}, \mathbf{a})$ and  $(\tilde{\boldsymbol{\psi}}, \tilde{\mathbf{a}})$
satisisfy $\mathbb{T}(\boldsymbol{\psi}, \mathbf{a})= (0,\mathbf{e}_1)$ and $\mathbb{T}(\tilde{\boldsymbol{\psi}}, \tilde{\mathbf{a}})= (0,\mathbf{e}_2)$ respectively, where the vectors 
$\mathbf{a}=[a_1, a_2]$ and $\tilde{\mathbf{a}}=[a_2, a_1]$. Then it follows that
\begin{equation}\label{eqn:Spsi_phi1}
\left[
\begin{array}{cc}
\langle \mathbb{S}  \boldsymbol{\psi}, \boldsymbol{\varphi} \rangle ] \\
\langle \mathbb{S}  \tilde{\boldsymbol{\psi}}, \boldsymbol{\varphi} \rangle ]
\end{array}
\right] = 
\left[
\begin{array}{cc}
\beta_0 + a_1 &  a_2  \\
a_2 & \beta_0 + a_1 
\end{array}
\right]  \cdot \int_I \boldsymbol{\varphi} \, dx 
= 
\left[
\begin{array}{cc}
\beta_0 + a_1 &  a_2  \\
a_2 & \beta_0 + a_1 
\end{array}
\right]  \cdot  
\left[
\begin{array}{cc}
\langle  \boldsymbol{\varphi}, \mathbf{e}_1 \rangle  \\
\langle  \boldsymbol{\varphi}, \mathbf{e}_2   \rangle
\end{array}
\right]  .
\end{equation}
On the other hand, 
\begin{equation}\label{eqn:Spsi_phi2}
\left[
\begin{array}{cc}
\langle \mathbb{S}  \boldsymbol{\psi}, \boldsymbol{\varphi} \rangle  \\
\langle \mathbb{S}  \tilde{\boldsymbol{\psi}}, \boldsymbol{\varphi} \rangle 
\end{array}
\right] = 
\left[
\begin{array}{cc}
\langle  \boldsymbol{\psi}, \mathbb{S} \boldsymbol{\varphi} \rangle  \\
\langle  \tilde{\boldsymbol{\psi}},\mathbb{S}  \boldsymbol{\varphi} \rangle 
\end{array}
\right] =
\left[
\begin{array}{cc}
\langle  \boldsymbol{\psi}, \mathbf{e}_1 \rangle  \\
\langle  \tilde{\boldsymbol{\psi}}, \mathbf{e}_1   \rangle
\end{array}
\right] = \mathbf{e}_1.
\end{equation}
From \eqref{eqn:Spsi_phi1} and \eqref{eqn:Spsi_phi2} and by repeating the calculation with $\tilde{\boldsymbol{\varphi}}$,
we obtain
\begin{equation*}
\left[
\begin{array}{cc}
\beta_0 + a_1 &  a_2  \\
a_2 & \beta_0 + a_1 
\end{array}
\right] \cdot \hat{\mathbb{Q}} = \mathbb{I}.
\end{equation*}
Therefore, by choosing $\beta_0$ satisfying $(\beta_0 + a_1)^2-a_2^2\neq0$, the proof is complete. \qed  \\

Recall that $\mathbb{L}_\pm = \mathbb{S} + \mathbb{S}_\pm^\infty$. It follows from \eqref{eq:decomp_op} - \eqref{eq:S_infty_norm2} that
$\norm {\Delta \mathbb{Q}_\pm } =\|\mathbb{Q}_{\pm} - \hat{\mathbb{Q}}\|=O(\tau(\kappa,k,\varepsilon))$ holds in $D_1$ and $D_2$, where
\begin{equation}\label{eq:tau}
\tau(\kappa,k,\varepsilon) = 
\left\{ \begin{array}{ll}
\varepsilon^2|\ln\varepsilon|+|\kappa|\varepsilon & \mbox{if} \; (\kappa, k) \in D_1;  \\
\varepsilon          & \mbox{if} \; (\kappa, k) \in D_2,
\end{array}
\right.
\end{equation}
Thus $\norm {\Delta  \mathbb{M}_\pm} = \| \mathbb{M}_{\pm} - \hat{\mathbb{M}}_{\pm} \| = O(\tau)$. We now are ready to apply the perturbation theory in
Section \ref{sec:perturb_mat} to investigate the roots of $\lambda_{j,\pm}(k; \kappa, \varepsilon)$ for the matrix $\mathbb{M}_{\pm}$.

Let
$$\alpha = \langle \mathbb{S}^{-1} \mathbf{e}_1, \mathbf{e}_1 \rangle = \langle\mathbb{S}^{-1} \mathbf{e}_2, \mathbf{e}_2 \rangle 
 \quad \mbox{and} \quad
 \tilde\alpha = \langle \mathbb{S}^{-1} \mathbf{e}_1, \mathbf{e}_2 \rangle = \langle\mathbb{S}^{-1} \mathbf{e}_2, \mathbf{e}_1 \rangle. $$
Define the function
$$\gamma(\kappa,k) = 2\beta_{\mathrm{e}}(\kappa,k,\varepsilon) + \dfrac{2}{\pi}\ln 2 - \dfrac{2}{\pi}\ln\varepsilon  - \beta_0. $$
First, using the formula \eqref{eq:eigenval1_M_hat} and recalling from \eqref{eq:beta} that $\beta=\beta_\mathrm{e}+\beta_\mathrm{i}-\beta_0$,
we obtain the eigenvalues of $\hat{\mathbb{M}}_\pm$:
\begin{eqnarray*}
\hat\lambda_{1,\pm}(k; \kappa,\varepsilon) &=& \varepsilon+\varepsilon (\beta \pm \tilde \beta + \beta_\mathrm{e}) (\alpha + \tilde\alpha)
=\varepsilon+\varepsilon (\beta_\mathrm{i} \pm \tilde \beta + 2\beta_\mathrm{e} - \beta_0) (\alpha + \tilde\alpha), \\
& =&\varepsilon+\varepsilon \left(\frac{1}{k\tan k} \pm \frac{1}{k\sin k} + \frac{2}{\pi}\varepsilon\ln\varepsilon + \varepsilon\gamma(\kappa,k)\right) (\alpha + \tilde\alpha), \\
\hat\lambda_{2,\pm}(k; \varepsilon) &=& \varepsilon+\varepsilon (\beta \pm \tilde \beta - \beta_\mathrm{e}) (\alpha - \tilde\alpha)
=\varepsilon+\varepsilon (\beta_\mathrm{i} \pm \tilde \beta - \beta_0) (\alpha - \tilde\alpha)\\
&=& \varepsilon+\varepsilon \left(\frac{1}{k\tan k} \pm \frac{1}{k\sin k} +\frac{2\ln 2}{\pi}- \beta_0 \right) (\alpha - \tilde\alpha).\end{eqnarray*}
%They are explicitly expressed by
%\begin{eqnarray*}
%\hat\lambda_{1,\pm}(k; \kappa,\varepsilon) &=& 
%=\varepsilon+\varepsilon \left(\frac{1}{k\tan k} \pm \frac{1}{k\sin k} + \frac{2}{\pi}\varepsilon\ln\varepsilon + \varepsilon\gamma(\kappa,k)\right) (\alpha + \tilde\alpha), \\
%\hat\lambda_{2,\pm}(k; \varepsilon) &=& \varepsilon+\varepsilon \left(\frac{1}{k\tan k} \pm \frac{1}{k\sin k} - \beta_0 \right) (\alpha - \tilde\alpha),
%\end{eqnarray*}
It is clear that, for real $k$,  $\hat\lambda_{1,\pm}(k)$ is complex-valued while $\hat\lambda_{2,\pm}(k)$ is real-valued.

The leading-order term of $\hat\lambda_{1,+}$ and $\hat\lambda_{2,+}$ 
is equal and is given by $\frac{1}{k\tan k}  + \frac{1 }{k\sin k}$.
In view of Lemma \ref{lem:perturb_eig_M}, we deduce that $\lambda_{1,+}$ and $\lambda_{2,+}$ share the same leading order, thus both attain
a simple root $m\pi$ for odd integers $m$.
Let us choose the disc $B_{\rho}(m\pi)$ with radius $\rho$ centered at $m\pi$ in the complex $k$-plane, with $\rho=O(1)$ as $\varepsilon\to0$. 
By extending the functions $\beta_\mathrm{e}(k)$,  $\beta_\mathrm{i}(k)$ and  $\tilde \beta $ to $B_{\rho}(m\pi)$,
one can show that the asymptotic expansions in $\varepsilon$ for the kernels $G_\varepsilon^\mathrm{e}$, $G_\varepsilon^\mathrm{i}$ and $\tilde G_\varepsilon^\mathrm{i}$ given in \eqref{Ge_exp}, \eqref{Gi_exp}, and \eqref{tGi_exp} hold in $B_{\rho}(m\pi)$.
From Rouche's theorem, $\lambda_{j,+}(k)$ ($j=1,2$) attains a simple root $k_{m}^{(j)}$
that is close to $m\pi$ if $\varepsilon$ is sufficiently small.

The asymptotics of $k_{m}^{(j)}$ can be obtained as follows. 
We apply the Taylor expansion for $\hat \lambda_{j,+}(k)$ at $k=m\pi$ to obtain the expansions for their roots:
% we let $\delta k:=k-m\pi$ and apply the Taylor expansion for $\hat \lambda_{j,+}(k)$ at $k=m\pi$:
% \begin{eqnarray*}
% \hat \lambda_{1,+}(k) &=& \varepsilon +\bigg[- \frac{1}{2m\pi} \cdot \delta k +O(\delta k^2) +
% \dfrac{2}{\pi} \varepsilon \ln \varepsilon + \varepsilon \gamma(m\pi, \kappa) + \varepsilon \cdot O(\delta k)  \bigg] \cdot (\alpha+\tilde\alpha); \\
% \hat \lambda_{2,+}(k) &=& \varepsilon +\bigg[- \frac{1}{2m\pi} \cdot \delta k +O(\delta k^2) + \varepsilon \left(\frac{2\ln 2}{\pi} - \beta_0\right) \bigg] \cdot (\alpha-\tilde\alpha).
% \end{eqnarray*}
% Thus the root of $\hat \lambda_{1,+}$ and $\hat \lambda_{2,+}(k)$ can be expanded as
\begin{eqnarray*}
\hat k_{m}^{(1)} &=& m \pi + 2m\pi \left [ \dfrac{2}{\pi}\varepsilon \ln\varepsilon + \left( \frac{1}{\alpha+\tilde\alpha} + \gamma(\kappa,m\pi) \right)\varepsilon \right] + \hat k_{m,h}^{(1)};  \\
\hat k_{m}^{(2)} &=& m \pi + 2m\pi \left(  \dfrac{1}{\alpha-\tilde\alpha} + \frac{2\ln 2}{\pi} - \beta_0 \right)\varepsilon  +
\hat k_{m,h}^{(2)},
\end{eqnarray*}
where $\hat k_{m,h}^{(1)}=O(\varepsilon^2\ln^2\varepsilon)$ and $\hat k_{m,h}^{(2)}=O(\varepsilon^2)$.
Furthermore, there holds
\begin{equation}\label{eq:Im_hat_k}
\Im\, \hat k_{m}^{(1)}  = O(\varepsilon) \quad \mbox{and} \quad \Im \, \hat k_{m}^{(2)} =0,
\end{equation}
since $\Im \, \gamma(\kappa,m\pi) \neq 0$ while $\hat\lambda_{2,+}(k)$ is a real-valued function for real $k$.

Now it follows from Lemma~\ref{lem:perturb_eig_M} that
$$  
\lambda_{j,+}(k) - \hat\lambda_{j,+}(k)  \;=\;  O(\tau)  \cdot   \hat\lambda_{j,+}(k)    + O(\varepsilon \tau)  \quad \mbox{in} \,\,\, B_{\delta_0}(m\pi).
 $$ 
%holds in $B_{\delta_0}(m\pi)$.
Hence there exists certain constant $C$ such that
$| \lambda_{j,+}(k) - \hat\lambda_{j,+}(k)|  < | \hat\lambda_{j,+}(k)| $
for all $k$ satisfying $|k-\hat k_{m}^{(j)}| = C |\varepsilon \tau|$. As such $k_{m}^{(j)}=\hat k_{m}^{(j)} + O(\varepsilon \tau)$.
Similarly, investigating the roots of $\hat \lambda_{j,-}(k; \kappa, \varepsilon)$ yields resonances
close to $m\pi$ with even integers~$m$.
In summary, we arrive at the following asymptotic expansion of resonances.

\begin{thm}\label{thm:asym_eig_perturb}
For each $\kappa \in (-b/2,b/2]$,
the scattering problem \eqref{eq:Helmholtz}--\eqref{eq:rad_cond2} admits two groups of complex-valued resonances
in the region $\mathbb{C}_\kappa$ given by
\begin{eqnarray}
k_m^{(1)} &=&  m \pi + 2m\pi \left [ \dfrac{2}{\pi}\varepsilon \ln\varepsilon + \left( \frac{1}{\alpha+\tilde\alpha} + \gamma(\kappa,m\pi) \right)\varepsilon \right] +k_{m,h}^{(1)}, \label{eq:km1} \\
k_m^{(2)} &=&  m \pi + 2m\pi \left(  \dfrac{1}{\alpha-\tilde\alpha} + \frac{2\ln 2}{\pi} - \beta_0 \right)\varepsilon  +
 k_{m,h}^{(2)} \label{eq:km2}
\end{eqnarray} 
for $m \varepsilon \ll 1$. In the above,
\begin{equation}\label{eq:kmh}
k_{m,h}^{(1)}=O(\varepsilon^2\ln^2\varepsilon), \quad \Re \, k_{m,h}^{(2)}=O(\varepsilon^2),\; \mbox{and} \; \Im \, k_{m,h}^{(2)}= O\Big(\varepsilon \tau(\kappa,m\pi,\varepsilon)\Big).
\end{equation}
\end{thm}

The imaginary part of $k_m^{(j)}$ plays the key role in the resonance effect  and we examine it closely in what follows.

\begin{prop}\label{prop:imag_res}
For each $\kappa \in (-b/2,b/2]$, let $k_m^{(1)}$ and $k_m^{(2)}$  be the two groups of resonances given in \eqref{eq:km1} and \eqref{eq:km2}.
Then their imaginary parts attain the following orders:
\begin{equation*}
\Im \, k_m^{(1)} = O(\varepsilon),  \quad
 \Im \, k_m^{(2)}=
 \left\{ \begin{array}{ll}
O(\kappa\varepsilon^2)& \mbox{if} \; (\kappa, m\pi) \in D_1;  \\
O(\varepsilon^2)          & \mbox{if} \; (\kappa, m\pi) \in D_2,
\end{array}
\right.
\end{equation*}
where $D_1$ and $D_2$ are defined in \eqref{eq:D1} and \eqref{eq:D2} respectively.
\end{prop}

\noindent\textbf{Proof.} First, a combination of \eqref{eq:Im_hat_k} and \eqref{eq:kmh} implies that
$$\Im \, k_m^{(1)} = O(\varepsilon) \quad \mbox{and} \quad 
  \Im \, k_m^{(2)}=\Im \, k_{m,h}^{(2)}= O(\varepsilon \tau). $$
It is clear that the order of $\Im \, k_m^{(1)}$ is optimal as $\varepsilon\to0$. 
If $(\kappa,m\pi) \in D_2$, from \eqref{eq:tau} it follows that $\Im \, k_m^{(2)}$ attains the order $O(\varepsilon^2)$.
On the other hand, if $(\kappa,m\pi) \in D_1$, one can only tell that
$\Im \, k_m^{(2)}=O(\varepsilon^3|\ln\varepsilon|+|\kappa|\varepsilon^2)$. In what follows, we apply the perturbation theory in Section \ref{sec:perturb_mat} recursively to 
obtain a sharper estimation and show that $\Im \, k_m^{(2)}=O(|\kappa|\varepsilon^2)$ if $(\kappa,m\pi) \in D_1$.

Without loss of generality, we focus on $k_m^{(2)}$ for odd $m$. Let $k_n$ $(n=0, 1, 2, \cdots)$ be a sequence of real numbers with $k_0=\hat k_{m}^{(2)}$ and $k_n  \in (0,b)$ for $n\ge1$ specified below.
Define matrice
\begin{equation}  \label{eq:matrix-Mn_hat}
\mathbb{M}_n:=\varepsilon(\mathbb{Q}_n\mathbb{B}_+ + \mathbb{I})
\quad \mbox{where} \;
\mathbb{Q}_n
=
\left[
\begin{array}{cc}
 \langle (\mathbb{L}_{n,0}^{-1}  \mathbf{e}_1,  \mathbf{e}_1 \rangle  &  \langle \mathbb{L}_{n,0}^{-1} \mathbf{e}_2,  \mathbf{e}_1\rangle \\
 \langle \mathbb{L}_{n,0}^{-1}  \mathbf{e}_1,  \mathbf{e}_2 \rangle  &  \langle \mathbb{L}_{n,0}^{-1}  \mathbf{e}_2,  \mathbf{e}_2 \rangle 
 \end{array}
\right],
\end{equation}
and $\mathbb{L}_{n,0}$ is the operator $\mathbb{L}_+(k,\kappa)$ for $k=k_n$ and $\kappa=0$. 
Similar to Lemma \ref{lem:Q_hat}, it can be shown that ${\mathbb{Q}}_n$ is real and symmetric. As such,
by letting 
$$ \alpha_n = \langle \mathbb{L}_{n,0}^{-1} \mathbf{e}_1, \mathbf{e}_1 \rangle = \langle\mathbb{L}_{n,0}^{-1} \mathbf{e}_2, \mathbf{e}_2 \rangle 
\quad \mbox{and} \quad 
\tilde\alpha_n = \langle \mathbb{L}_{n,0}^{-1} \mathbf{e}_1, \mathbf{e}_2 \rangle = \langle\mathbb{L}_{n,0}^{-1} \mathbf{e}_2, \mathbf{e}_1 \rangle, $$
the second eigenvalue of $\mathbb{M}_{n}$ reads
 \begin{equation}
 \lambda_{n,2}(k; \kappa,\varepsilon) = \varepsilon+\varepsilon (\beta + \tilde \beta - \beta_\mathrm{e}) (\alpha_n - \tilde\alpha_n),
 \label{eq:eigenval2_M_bar}
 \end{equation}
which attains real root $k_{n+1}$ near $k_m^{(2)}$. A close examination using the estimation in Lemma \ref{lem:hot} yields
\begin{eqnarray}
&& |k_1-k_0|=O(\varepsilon \cdot (\alpha_{0} - \alpha))=O(\varepsilon^3|\ln\varepsilon|),  \label{eq:kn1} \\
&& |k_{n+1}-k_n|=O(\varepsilon \cdot |\alpha_{n} - \alpha_{n-1}|) = O(\varepsilon^3|\ln\varepsilon|) \cdot |k_n-k_{n-1}|. \label{eq:kn2}
\end{eqnarray}

Recall that $\mathbb{Q}_{+}(k):= \mathbb{Q}_{+}(k, \kappa, \varepsilon)$ in (\ref{eq:matrix_QB}) and $\lambda_{2,+} := \lambda_{2,+} (k; \kappa,\varepsilon)$ in (\ref{eq:lambdas}). 
Using the relation
$$ \Delta  \mathbb{Q}_n := \mathbb{Q}_{+}(k) - \mathbb{Q}_n = \Big(\mathbb{Q}_{+}(k) - \mathbb{Q}_{+}(k_n)\Big) +\Big(\mathbb{Q}_{+}(k_n) -  \mathbb{Q}_n \Big), $$
it follows from \eqref{eq:S_infty_norm2} that $\norm{\Delta  {\mathbb{Q}}_n }=O(\varepsilon^2|\ln\varepsilon|) \cdot (k-k_n) +O(\kappa\varepsilon)$.
An application of Lemma \ref{lem:perturb_eig_M} yields
$$\lambda_{2,+} (k; \kappa,\varepsilon)  = (1+O(\tau_n) ) \cdot \lambda_{n,2} (k; \kappa,\varepsilon)  + O(\varepsilon \tau_n), 
\quad \mbox{where} \; \tau_n:=\varepsilon^2|\ln\varepsilon| \cdot (k-k_n)  + \kappa\varepsilon. $$
This is can be rewritten as
$$  \lambda_{2,+}(k) - \lambda_{n,2}(k)  = 
O(\tau_n)  \cdot  \lambda_{n,2}(k)     +  O(\varepsilon^3|\ln\varepsilon|) \cdot \big[(k-k_{n+1})  + (k_{n+1}-k_n)\big]  + O(\kappa\varepsilon^2), $$
which implies that
$| \lambda_{2,+}(k) - \lambda_{n,2}(k)|  < |\lambda_{n,2}(k)| $
for all $k$ satisfying 
\begin{equation}\label{eq:dis_k_kn1}
|k-k_{n+1}| = C (\varepsilon^3|\ln\varepsilon| \cdot |k_{n+1}-k_n| + \kappa\varepsilon^2)
\end{equation}
for certain constant $C$. As such $k_m^{(2)}$, the root of $ \lambda_{2,+}(k)$, sits in the disk with the above radius by Rouche's theorem.
In light of \eqref{eq:kn1} - \eqref{eq:dis_k_kn1}, for any given nonzero $\kappa$, there exists sufficiently large integer $N$ such that  $|k_m^{(2)}-k_{N}| \lesssim |\kappa|\varepsilon^2$
and there holds $ \Im \, k_m^{(2)}=O(\kappa\varepsilon^2)$.

Finally, if $\kappa=0$, $\mathbb{Q}_\pm$ is real and symmetric for real-valued $k$, so does $\lambda_{2,+}(k)$. Thus $k_m^{(2)}$ is real. This completes the proof of the proposition.

\section{Fano resonance for the grating structure}\label{sec:anomaly}
\subsection{Asymptotics of the solution to scattering problem}\label{subsec:scattering}
We solve the system $\mathbb{T}\boldsymbol{\varphi}=\varepsilon^{-1}\mathbf{f}$ with $\mathbb{T}$, $\boldsymbol{\varphi}$ and $\mathbf{f}$
given in \eqref{eq:system_BIE}. To this end, we apply the decomposition $\mathbf{f}=\mathbf{f}_\mathrm{even}+\mathbf{f}_\mathrm{odd}$, 
with ${\mathbf{f}}_\mathrm{even} = [f^-, \, f^+, \, f^-, \, f^+ ]^T$ and ${\mathbf{f}}_\mathrm{odd} = [f^-, \, f^+, \, -f^-, \, -f^+ ]^T$,
and solve the two subsystems
$\mathbb{T}\boldsymbol{\varphi}_\mathrm{even}=\varepsilon^{-1}{\mathbf{f}}_\mathrm{even}$ and
$\mathbb{T}\boldsymbol{\varphi}_\mathrm{odd}=\varepsilon^{-1}{\mathbf{f}}_\mathrm{odd}$ separately.
From the structure of the operator $\mathbb{T}$, one can express
\begin{equation}\label{eq:phi_eo}
\boldsymbol{\varphi}_\mathrm{even} = [ \boldsymbol{\varphi}_{+}, \boldsymbol{\varphi}_{+}]^T
\quad \mbox{and} \quad
\boldsymbol{\varphi}_\mathrm{odd} = [ \boldsymbol{\varphi}_{-}, -\boldsymbol{\varphi}_{-}]^T,
\end{equation}
then these two subsystems are equivalent to the two smaller systems
$$  \mathbb{T}_+\boldsymbol{\varphi}_{+} =  \varepsilon^{-1}\tilde{\mathbf{f}} \quad \mbox{and} \quad  \mathbb{T}_-\boldsymbol{\varphi}_{-} =  \varepsilon^{-1}\tilde{\mathbf{f}},   $$
where $\mathbb{T}_+=\hat{\mathbb{T}}+\tilde{\mathbb{T}}$ and $\mathbb{T}_-=\hat{\mathbb{T}}-\tilde{\mathbb{T}}$
are defined in \eqref{eq:T_pm}, and $\tilde{\mathbf{f}} = [f^-, \, f^+]^T$.

The expansion of their solutions is given in the following lemma.
We distinguish the two cases: $(\kappa, k) \in D_1$ and $(\kappa, k) \in D_2$. We define 
\begin{equation*}
 \eta_{\pm}(k) =
 \left\{ \begin{array}{ll}
\big\langle [\mathbb{L}_{\pm}(0,k)]^{-1} \mathbf{e}_1, \mathbf{e}_1+\mathbf{e}_2\big\rangle & \mbox{if} \; (\kappa, k) \in D_1;  \\
\langle \mathbb{S}^{-1} \mathbf{e}_1, \mathbf{e}_1 + \mathbf{e}_2 \rangle        & \mbox{if} \; (\kappa, k) \in D_2,
\end{array}
\right.
\end{equation*}
and
\begin{equation*}
\delta(\varepsilon)=
 \left\{ \begin{array}{ll}
|\kappa|\varepsilon & \mbox{if} \; (\kappa, k) \in D_1;  \\
\varepsilon         & \mbox{if} \; (\kappa, k) \in D_2.
\end{array}
\right.
\end{equation*}

\begin{lem} \label{lem:phi}
The following asymptotic expansion holds for the solutions $\boldsymbol{\varphi}_+$ and  $\boldsymbol{\varphi}_-$ in $V_1\times V_1$:
\begin{equation}\label{eq:phi_1}
\left[
\begin{array}{llll}
\langle  \boldsymbol{\varphi}_\pm, \mathbf{e}_1 \rangle  \\
\langle  \boldsymbol{\varphi}_\pm, \mathbf{e}_2 \rangle
\end{array}
\right]  \nonumber
= - \dfrac{1}{\lambda_{1,\pm}} \left( \eta_{\pm}+O(\delta) \right)\cdot \mathbf{v}  
+ \dfrac{1}{\lambda_{2,\pm}} \cdot \mathbf{w}_{\pm} \quad \mbox{if} \; (\kappa,k) \in D_j,
\end{equation} 
where $\mathbf{v}= [1, \, 1]^T$ and $\mathbf{w}_\pm=[w_{1,\pm}, \, w_{2,\pm}]$ with $w_{j,\pm}=O(\delta)$ ($j=1, 2$).
\end{lem}

\noindent\textbf{Proof.} Using the decomposition $\mathbb{T}_\pm=\mathbb{P}_\pm+\mathbb{L}_\pm$ and applying $\mathbb{L}_\pm^{-1}$ on both sides,
the equation $\mathbb{T}_\pm \boldsymbol{\varphi}_{+} = \varepsilon^{-1}\tilde{\mathbf{f}}$ can be expressed as
\begin{equation}\label{eq:sol_phi}
\mathbb{L}_\pm^{-1} \;\mathbb{P}_\pm \; \boldsymbol{\varphi}_\pm + \boldsymbol{\varphi}_\pm \;=\; \varepsilon^{-1} \mathbb{L}_\pm^{-1}  \tilde{\mathbf{f}}. 
\end{equation}
A calculation analogous to that in Section~\ref{sec:res_cond} leads to
\begin{equation}\label{eq:linear_system}
\mathbb{M}_{\pm} \left[
 \begin{array}{llll}
\langle  \boldsymbol{\varphi}_\pm, \mathbf{e}_1 \rangle  \\
\langle  \boldsymbol{\varphi}_\pm, \mathbf{e}_2 \rangle
\end{array}
\right] \;=\;
\left[
\begin{array}{llll}
 \langle \mathbb{L}_\pm^{-1} \tilde{\mathbf{f}}, \mathbf{e}_1 \rangle   \\
 \langle \mathbb{L}_\pm^{-1} \tilde{\mathbf{f}}, \mathbf{e}_2 \rangle 
\end{array}
\right].
\end{equation}

Recall that the matrix $\mathbb{M}_{\pm}$ has eigenvalues 
$\lambda_{1,\pm}$ and $\lambda_{2,\pm}$, which are associated with the eigenvectors $v_{1,\pm}$ and $v_{2,\pm}$.
We define the matrix $\hat{\mathbb{M}}_\pm$ in the form of \eqref{eq:matrix-M_hat}, where the symmetric matrix $\hat{\mathbb{Q}}$
is now defined by
\begin{equation*}
[\hat{\mathbb{Q}}]_{ij} = \big\langle [\mathbb{L}_{\pm}(0, \hat k)]^{-1} \mathbf{e}_i, \mathbf{e}_j \big\rangle 
\quad \mbox{and} \quad
[\hat{\mathbb{Q}}]_{ij} =  \langle \mathbb{S}^{-1}  \mathbf{e}_i,  \mathbf{e}_j \rangle
\end{equation*}
for $(\kappa,k) \in D_1$ and $D_2$ respectively. It follows from the relation \eqref{eq:re_D1} that $\| \mathbb{Q} - \hat{\mathbb{Q}} \| = O(\delta)$.
In light of Lemma \ref{lem:perturb_eig_M}, one can write $v_{1,\pm}=\hat v_{1} + \Delta v_{1,\pm} $ and 
$v_{2,\pm}=\hat v_{2} + \Delta v_{2,\pm}$, where $\hat v_{1}$ and $\hat v_{2}$ are eigenvectors of
$\hat{\mathbb{M}}_\pm$ given explicitly in \eqref{eq:eigenvec_M_hat}. A direct calculation yields
\begin{equation}\label{eq:M_inv}
\mathbb{M}_{\pm}^{-1} = \dfrac{1}{2\lambda_{1,\pm}}
\left(\left[
\begin{array}{cc}
1 & 1 \\
1 & 1 
\end{array}
\right] + \Delta \mathbb{V}_{1,\pm} \right)
+\dfrac{1}{2\lambda_{2,\pm}}
\left(\left[
\begin{array}{cc}
1 & -1 \\
-1 & 1 
\end{array}
\right] + \Delta \mathbb{V}_{2,\pm} \right),
\end{equation}
where $\Delta\mathbb{V}_{j,\pm}$ is a $2\times2$ matrix with $\norm{\Delta\mathbb{V}_{j,\pm}} = \norm{\Delta v_{j,\pm}} = O(\delta)$ ($j=1,2$).
On the other hand, there holds
$$
\tilde{\mathbf{f}}= - (1+O(\kappa\varepsilon)) (\mathbf{e}_1+\mathbf{e}_2) ,  $$
and consequently,
\begin{equation}\label{eq:L_inv_f}
\left[
\begin{array}{llll}
 \langle \mathbb{L}_\pm^{-1} \tilde{\mathbf{f}}, \mathbf{e}_1 \rangle   \\
 \langle \mathbb{L}_\pm^{-1} \tilde{\mathbf{f}}, \mathbf{e}_2 \rangle 
\end{array}
\right] = -\left(\eta_{\pm}+O(\kappa\varepsilon)\right)
\left(\left[
\begin{array}{llll}
 1\\
 1
\end{array}
\right] + \Delta\tilde{\mathbf{f}} \right),
\end{equation}
where $\|\Delta\tilde{\mathbf{f}}\|=O(\delta)$.
The proof is completed by inserting \eqref{eq:M_inv} and \eqref{eq:L_inv_f} into \eqref{eq:linear_system}.  \qed \\

\begin{figure}
\begin{center}
\vspace*{-0.5cm}
\includegraphics[height=5cm]{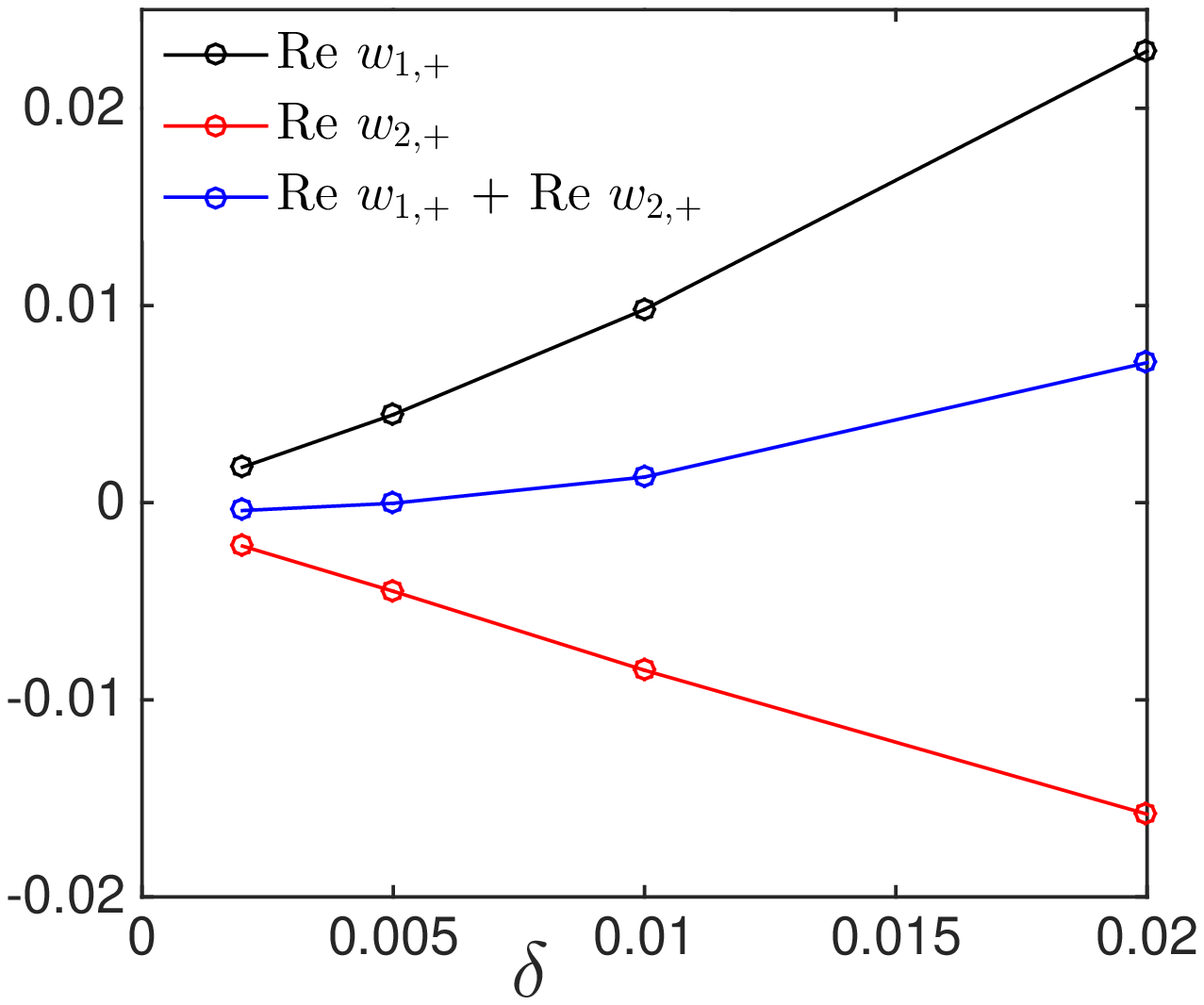}
\includegraphics[height=5cm]{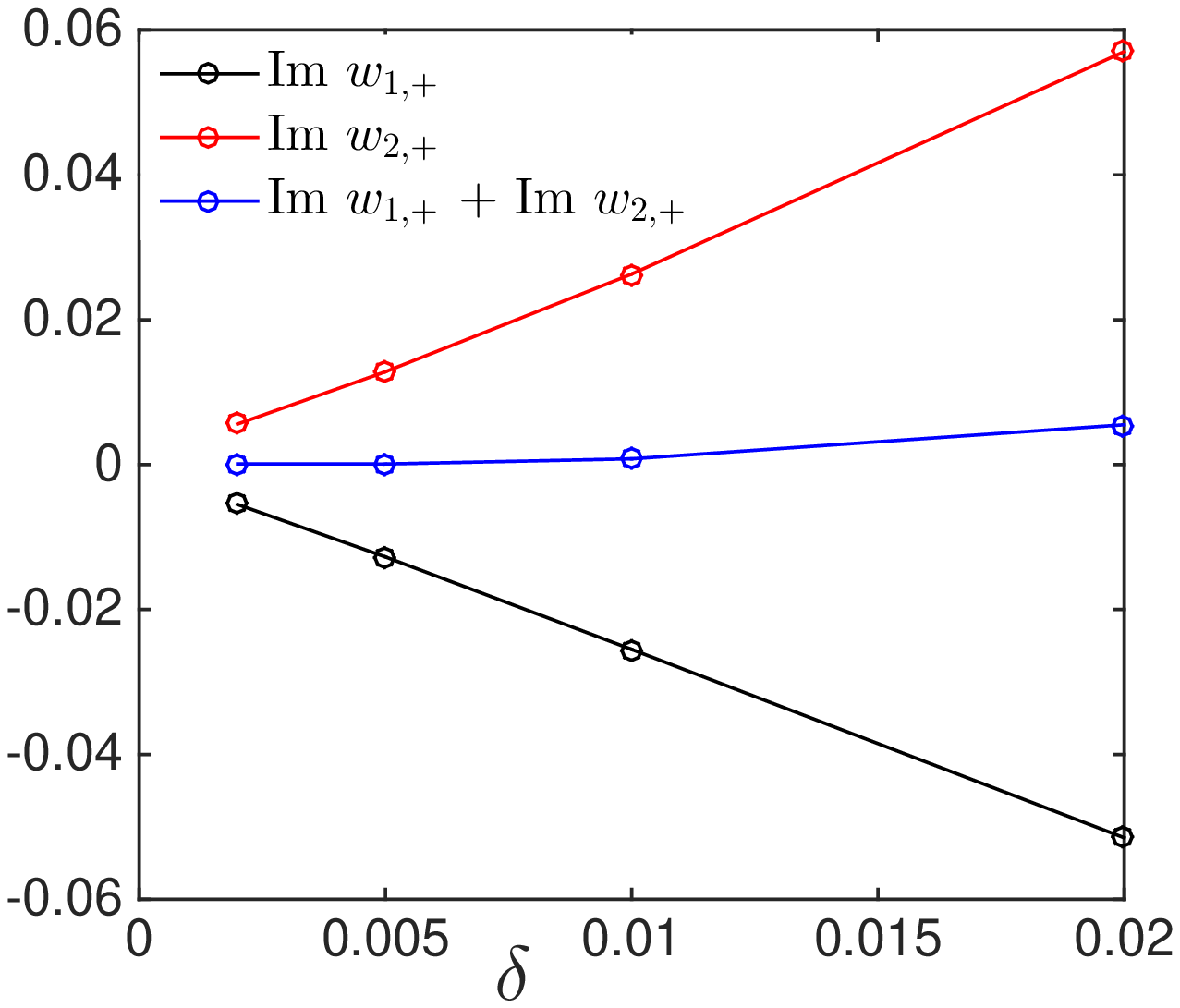} 
\caption{Real and imaginary parts of $w_{1,+}$,  $w_{2,+}$, and $w_{1,+}+w_{2,+}$  at the frequency $k=3.12$ for various $\delta$ values.
}\label{fig:wri}
\end{center}
\end{figure}

According to the lemma above, $w_{1,\pm}$ and $w_{2,\pm}$ have the same order $O(\delta)$. We would like to point out that their phase difference 
$\arg w_{1,\pm}-\arg w_{2,\pm}$ is about $\pi$ near the resonance frequencies $Re \, k_m^{(2)}$, and
there holds that
\begin{equation}\label{eq:sum_w1w2}
w_{1,\pm}+w_{2,\pm}=o(\delta).
\end{equation}
This has been observed through a large number of numerical calculations, although a rigorous proof is lacking here. Here we demonstrate such a relation in Figure \ref{fig:wri} for a structure with period $d=1.5$.

Let us now consider the reflected and transmitted wave above and below the grating. 
From the Green's formula, the total field above the grating is
$$ u_\varepsilon(x) = \int_{\Gamma^+_{1,\varepsilon} \cup \Gamma^-_{1,\varepsilon}} g_1(x,y) \partial_\nu u_\varepsilon(y) ds_y + u^\mathrm{inc}+ u^\mathrm{refl} \quad \mbox{in} \; \Omega_1. $$
For $x$ sufficiently far away from the grating, a change of variable to the scaled interval $I$ and the Taylor expansion leads to
$$
u_\varepsilon(x) = \varepsilon \Big( g_1\big(x,(-\varepsilon,1)\big)+O(\varepsilon) \Big) \cdot  \langle \varphi_1^-, 1 \rangle 
+ \varepsilon \Big( g_1\big(x,(\varepsilon,1)\big)+O(\varepsilon) \Big) \cdot  \langle \varphi_1^+, 1 \rangle  + u^\mathrm{inc}+ u^\mathrm{refl}.
$$
In view of the decomposition \eqref{eq:phi_eo} and the asymptotic expansion in Lemma~\ref{lem:phi}, for $(\kappa, k) \in D_j$, we obtain
\begin{equation}\label{eq:u_omega1}
u_\varepsilon(x) =  \varepsilon \hat r_- \cdot   \Big( g_1\big(x,(-\varepsilon,1)\big)+O(\varepsilon) \Big) 
+ \varepsilon \hat r_+ \cdot   \Big( g_1\big(x,(\varepsilon,1)\big)+O(\varepsilon) \Big)    + u^\mathrm{inc}+ u^\mathrm{refl},
\end{equation}
where the reflection coefficients
\begin{eqnarray}
\hat r_- (\kappa,k,\varepsilon) &:=&  \langle \varphi_1^-, 1 \rangle =  -\dfrac{1}{\lambda_{1,+}} \Big( \eta_{+} +O(\delta) \Big)  
- \dfrac{1}{\lambda_{1,-}} \Big( \eta_{-} +O(\delta) \Big)
 + \dfrac{w_{1,+}}{\lambda_{2,+}} + \dfrac{w_{1,-}}{\lambda_{2,-}},  \label{eq:r-} \\
\hat r_+ (\kappa,k,\varepsilon) &:=& \langle \varphi_1^+, 1 \rangle =  -\dfrac{1}{\lambda_{1,+}} \Big( \eta_{+} +O(\delta) \Big)  
- \dfrac{1}{\lambda_{1,-}} \Big( \eta_{-} +O(\delta) \Big)
 + \dfrac{w_{2,+}}{\lambda_{2,+}} + \dfrac{w_{2,-}}{\lambda_{2,-}}. \label{eq:r+} 
\end{eqnarray}
Similarly, for $x\in\Omega_2$,
\begin{equation} \label{eq:u_omega2}
u_\varepsilon(x) = \varepsilon \hat t_- \cdot   \Big( g_2\big(x,(-\varepsilon,0)\big)+O(\varepsilon) \Big) 
+ \varepsilon \hat t_+ \cdot   \Big( g_2\big(x,(\varepsilon,0)\big)+O(\varepsilon) \Big),
\end{equation}
where the transmission coefficients
\begin{eqnarray}
\hat t_- (\kappa,k,\varepsilon) &:=& \langle \varphi_2^-, 1 \rangle = -\dfrac{1}{\lambda_{1,+}} \Big( \eta_{+} +O(\delta) \Big)  
+ \dfrac{1}{\lambda_{1,-}} \Big( \eta_{-} +O(\delta) \Big)
 + \dfrac{w_{1,+}}{\lambda_{2,+}} - \dfrac{w_{1,-}}{\lambda_{2,-}},   \label{eq:t-}  \\
\hat t_+  (\kappa,k,\varepsilon) &:=& \langle \varphi_2^+, 1 \rangle =  -\dfrac{1}{\lambda_{1,+}} \Big( \eta_{+} +O(\delta) \Big)  
+ \dfrac{1}{\lambda_{1,-}} \Big( \eta_{-} +O(\delta) \Big)
 + \dfrac{w_{2,+}}{\lambda_{2,+}} - \dfrac{w_{2,-}}{\lambda_{2,-}}.   \label{eq:t+} 
\end{eqnarray}

\begin{rmk}\label{rmk:field}
From formula \eqref{eq:u_omega1} and \eqref{eq:u_omega2}, we observe that the diffracted field generated by the slits can be approximately viewed as the field generated
by equivalent point charges sitting on the slit apertures. The strength of the charges is determined by the coefficients $\hat r_\pm$ and $\hat t_\pm$.
\end{rmk}

\subsection{Fano-type transmission anomalies}\label{subsec:anomaly}
By decomposing the Green functions $g_1(x,y)$ and $g_2(x,y)$ into the propagating and exponentially decaying parts,
we obtain the expansion of the reflected and transmitted fields via the formulas \eqref{eq:u_omega1} and \eqref{eq:u_omega2}:
$$ u_\varepsilon^r(x) = \sum_{n\in\mathbb{Z}_1}  r_n e^{i\kappa_n x_1 + i\zeta_n(x_2-1)}, \quad
 u_\varepsilon^t(x) = \sum_{n\in\mathbb{Z}_1}  t_n e^{i\kappa_n x_1 - i\zeta_nx_2}, $$
where the reflection and transmission coefficients of diffraction order $n$ are given by
 \begin{eqnarray}\label{eq:rn_tn}
 r_n(\kappa,k,\varepsilon)  &=& 
 \left\{ \begin{array}{lll}
 1 - \dfrac{i\varepsilon}{d\,\zeta_0} \big( e^{i\kappa \varepsilon} \hat r_- + e^{-i\kappa \varepsilon} \hat r_+ \big) \cdot \left(1+O(\varepsilon)\right)  &  \quad n = 0;  \\
  -\dfrac{i\varepsilon}{d\,\zeta_n} \big( e^{i\kappa_n \varepsilon} \hat r_- + e^{-i\kappa_n \varepsilon} \hat r_+ \big) \cdot \left(1+O(\varepsilon)\right)  & \quad n \in \mathbb{Z}_1 \backslash \{ 0 \},
 \end{array}
 \right.   \label{eq:rn} \\
 t_n(\kappa,k,\varepsilon)  &=& -\dfrac{i\varepsilon}{d\,\zeta_n} \big( e^{i\kappa_n \varepsilon} \hat t_- + e^{-i\kappa_n \varepsilon} \hat t_+  \big) \cdot \left(1+O(\varepsilon)\right). 
 \quad\quad  n \in \mathbb{Z}_1.   \label{eq:tn}
 \end{eqnarray}
Assuming that the incident power is one, then the reflection and transmission power is given by
\begin{equation}\label{eq:RT}
|R(\kappa,k,\varepsilon)|^2 = \sum_{n\in\mathbb{Z}_1} \frac{\zeta_n}{\zeta_0}|r_n|^2   
\quad \mbox{and} \quad
|T(\kappa,k,\varepsilon)|^2 = \sum_{n\in\mathbb{Z}_1} \frac{\zeta_n}{\zeta_0}|t_n|^2
\end{equation}
respectively.

Now let us study the transmission intensity of the scattering problem in the vicinity of the resonance frequency $k_*:=\Re \, k_m^{(2)}$.
We focus the discussion for odd integer $m$ only since the calculations are parallel for even integer $m$. There holds 
$\varepsilon/\lambda_{1,-}=O(\varepsilon)$ and $\varepsilon/\lambda_{2,-}=O(\varepsilon)$ in the $O(\varepsilon\delta)$ neighborhood of $k_*$.
Hence, from \eqref{eq:r-}-\eqref{eq:r+} and  \eqref{eq:t-}-\eqref{eq:t+}, it follows that
\begin{eqnarray}
\varepsilon \hat r_- &=& -\dfrac{\varepsilon \eta_{+}}{\lambda_{1,+}}  + 
\dfrac{\varepsilon w_{1,+}}{\lambda_{2,+}} + O(\varepsilon),
\quad
\varepsilon \hat r_+ = -\dfrac{\varepsilon \eta_{+}}{\lambda_{1,+}} + 
\dfrac{\varepsilon w_{2,+}}{\lambda_{2,+}} + O(\varepsilon); \label{r_hat}  \\
\varepsilon \hat t_- &=& -\dfrac{\varepsilon \eta_{+}}{\lambda_{1,+}} + 
\dfrac{\varepsilon w_{1,+}}{\lambda_{2,+}} + O(\varepsilon),
\quad
\varepsilon \hat t_+ = -\dfrac{\varepsilon \eta_{+}}{\lambda_{1,+}} + 
\dfrac{\varepsilon w_{2,+}}{\lambda_{2,+}} + O(\varepsilon).  \label{t_hat}
\end{eqnarray}
In view of the asymptotic expansion for the resonances frequencies in Theorem~\ref{thm:asym_eig_perturb}, we have 
$\dfrac{\varepsilon \eta_{+}}{\lambda_{1,+}}=c_0 (1+O(\varepsilon))$ in the $O(\varepsilon\delta)$ neighborhood of $k_*$,
where the constant $c_0=O(1/|\ln\varepsilon|)$. On the other hand, a Taylor expansion at $k_*$ yields
\begin{equation}\label{eq:lambda2_exp}
\lambda_{2,+}(k) = c_1 (k- k_*) + i \, c_2 \, \hat\delta + O(\varepsilon) \cdot( k-k_*) = \left(c_1 (k-k_*) + i \, c_2 \, \hat\delta\right) \cdot (1+O(\varepsilon)),
\end{equation}
where $c_1$ and $c_2$ are constants. As such, there holds $\dfrac{\varepsilon w_{1,+}}{\lambda_{2,+}}=O(1)$ and $\dfrac{\varepsilon w_{2,+}}{\lambda_{2,+}}=O(1)$ in the $O(\varepsilon\delta)$ neighborhood of $k_*$ when the resonance is excited.
In addition, their phases change rapidly due to the presence of the factor $(k- k_*)$ in the expansion \eqref{eq:lambda2_exp}.

\begin{figure}
\begin{center}
\vspace*{-0.5cm}
\includegraphics[height=4.5cm]{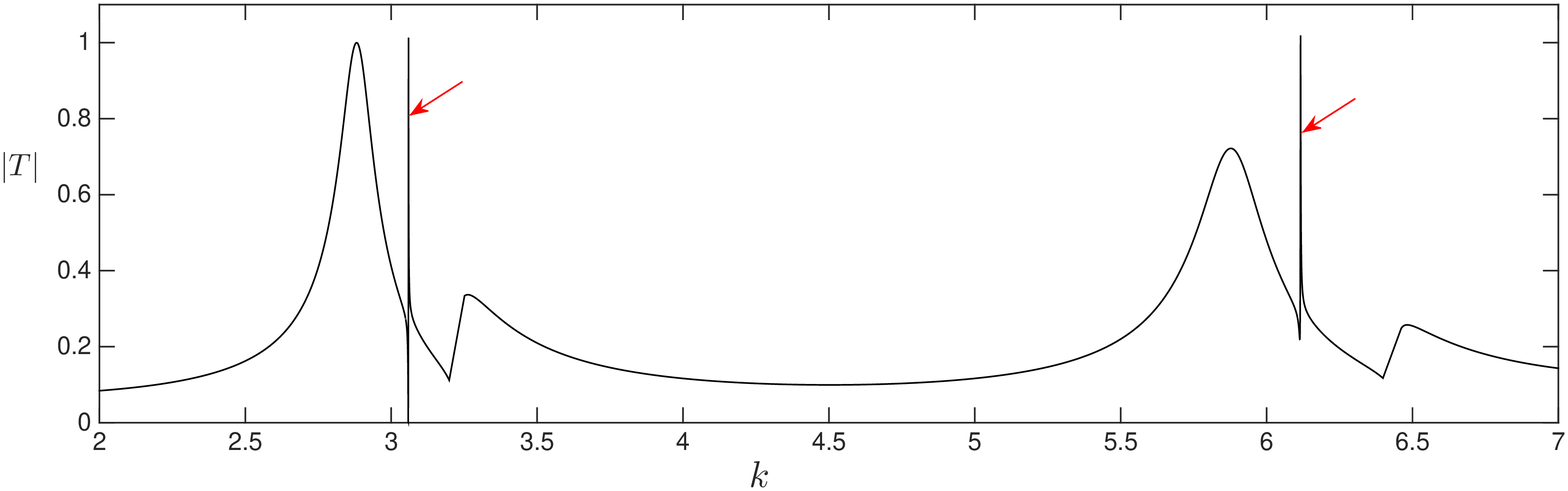}
\includegraphics[height=4cm]{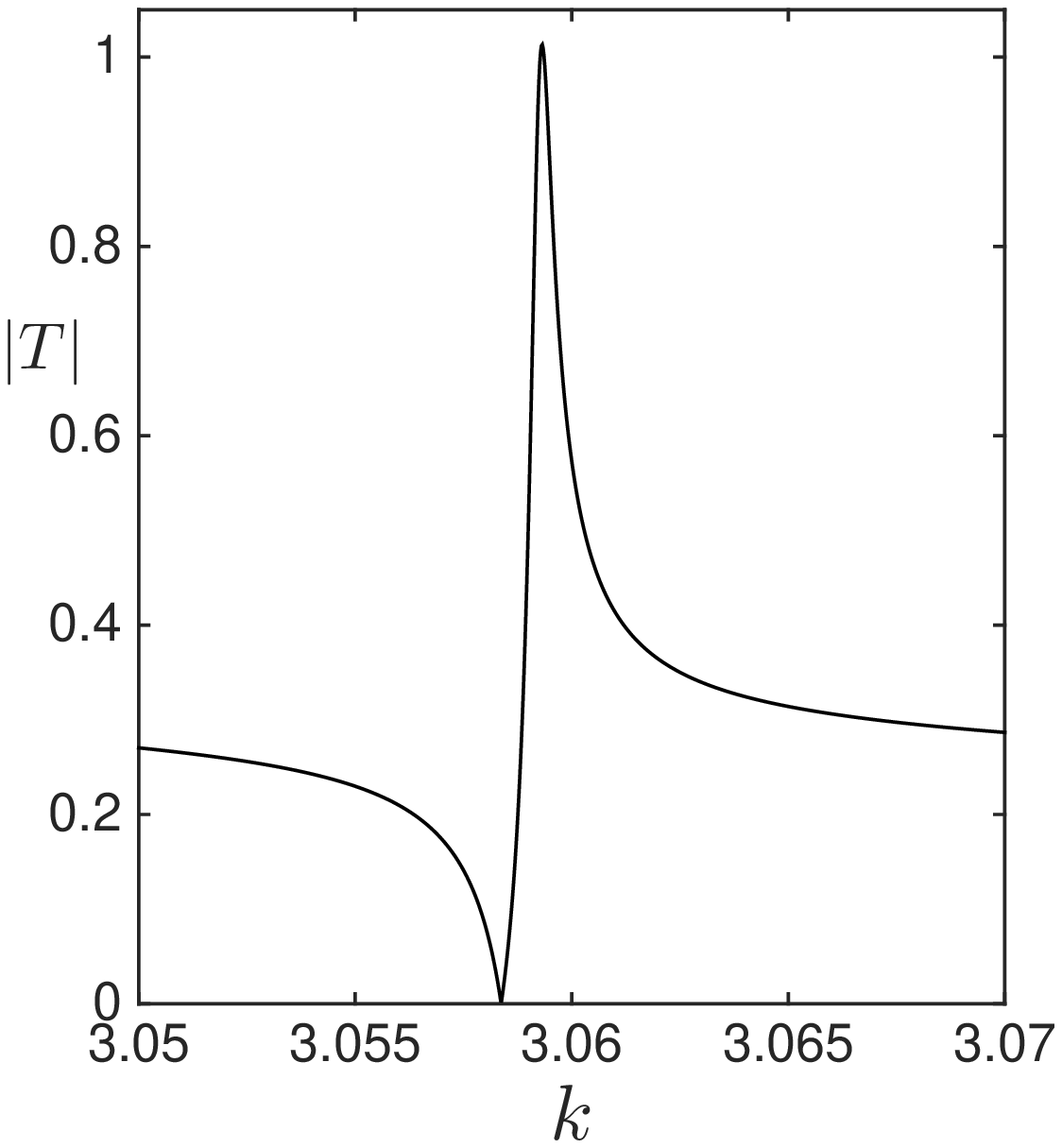} \quad\quad
\includegraphics[height=4cm]{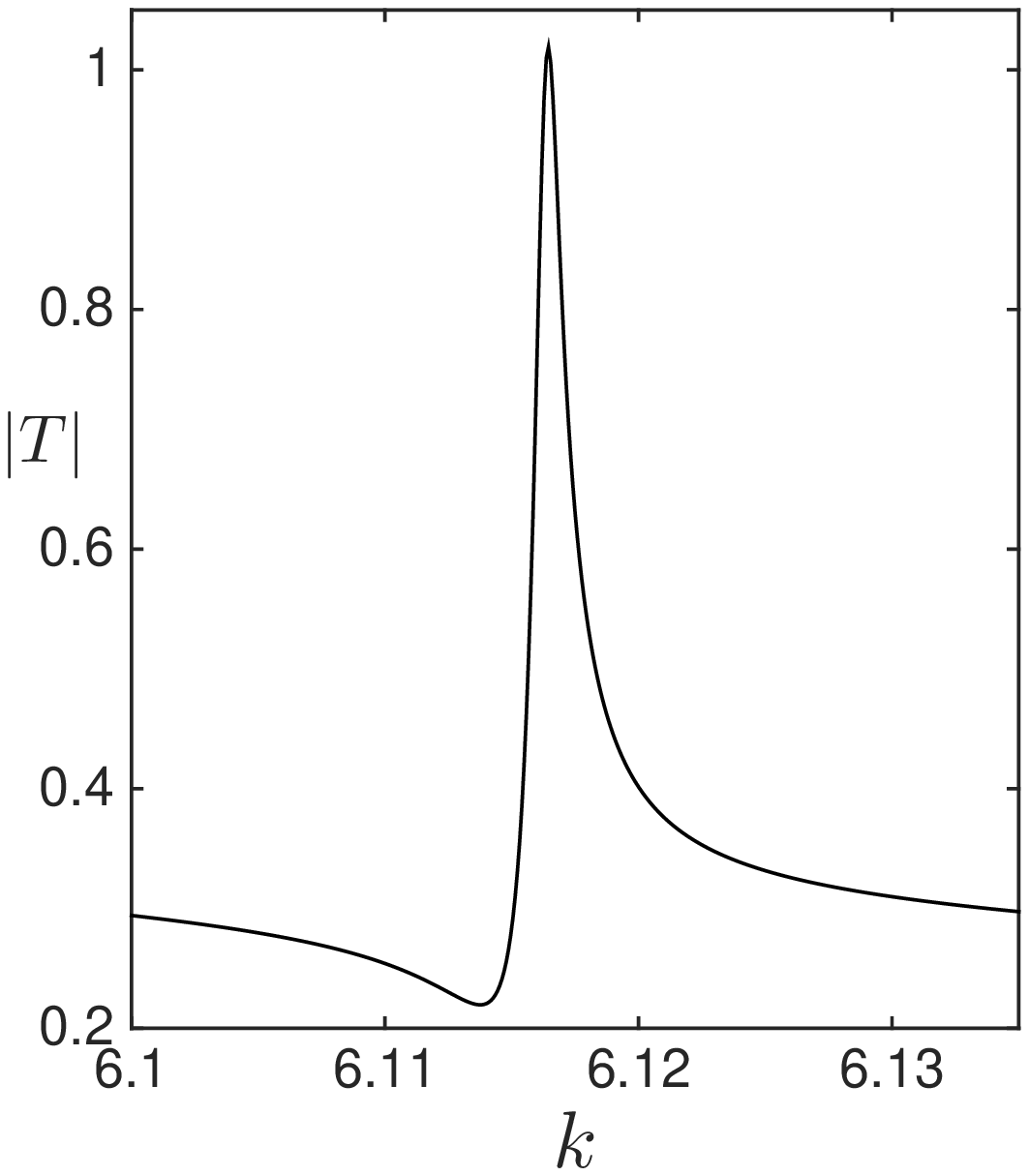}
\caption{Top: Transmission $|T|$ for $k \in (2,7)$ when $d=1.3$, $\varepsilon=0.02$. The incident angle $\theta=\pi/6$. Fano resonance (highlighted with red arrows)
occurs near $k=3.06$ and $6.12$.  Bottom: zoomed view of the two Fano-type transmission anomalies.
Two Fano resonances occurs for $(\kappa, k) \in D_1$ and $D_2$ respectively.
}\label{fig:trans_fano1}
\end{center}
\end{figure}

\begin{figure}
\begin{center}
\includegraphics[height=4.8cm]{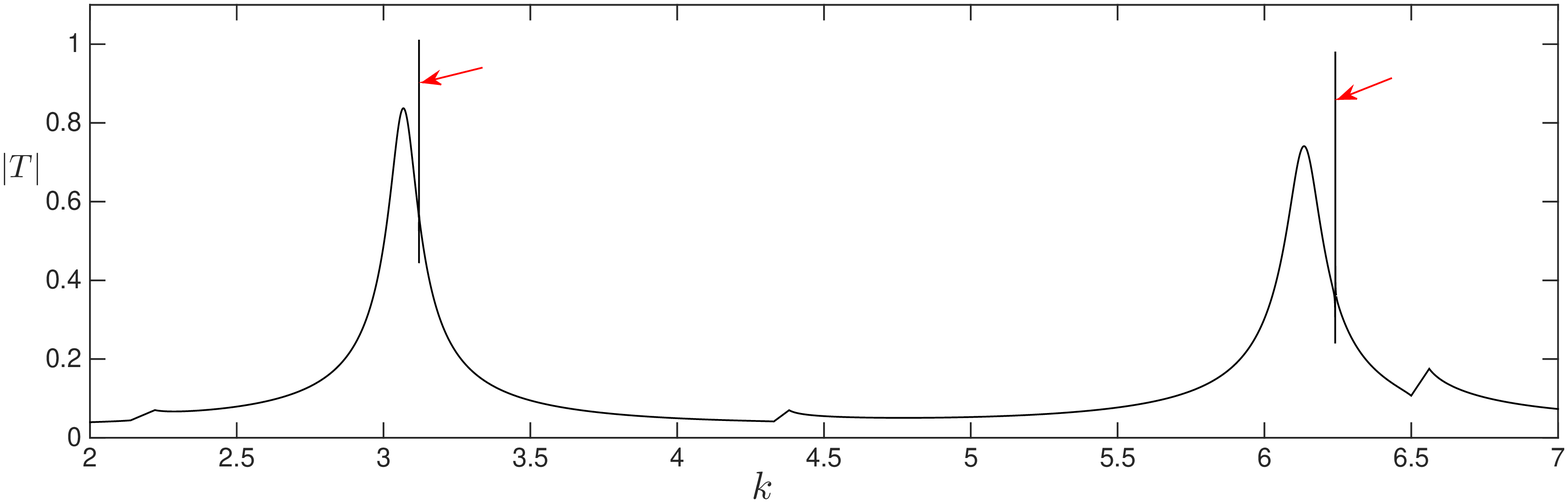}
\includegraphics[height=4cm]{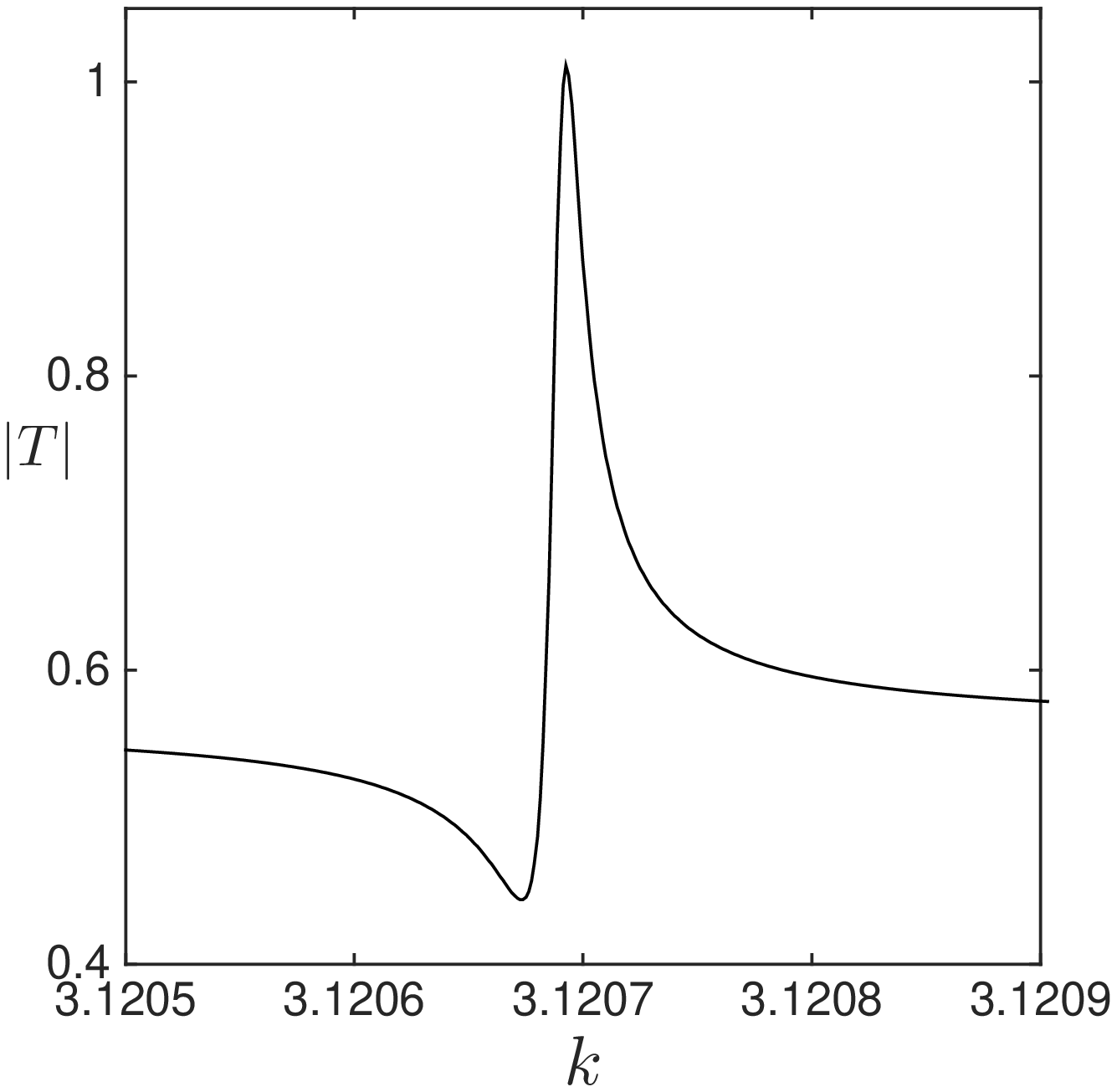}  \quad\quad
\includegraphics[height=4cm]{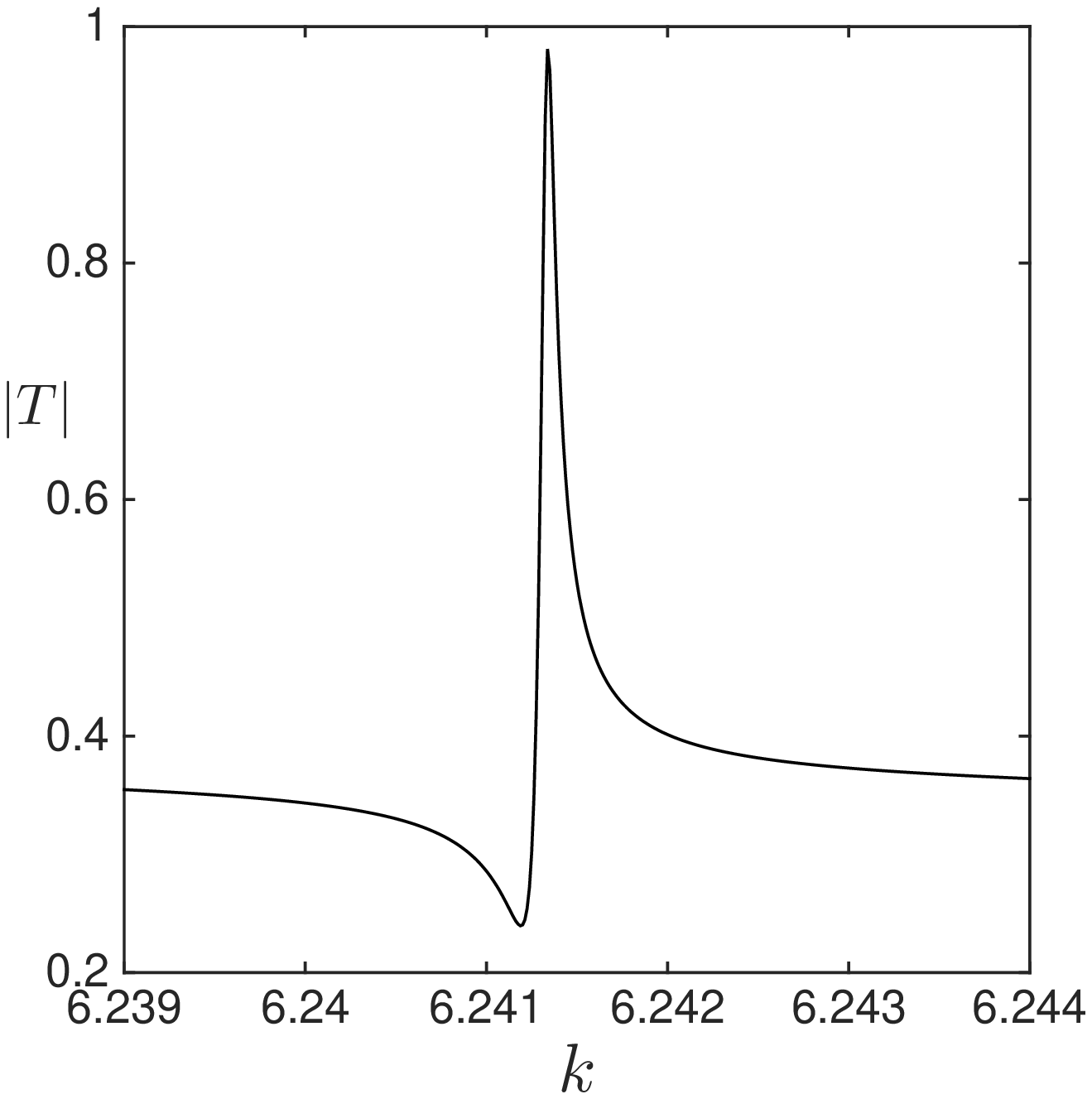}
\caption{Top: Transmission $|T|$ for $k \in (2,7)$ when $d=1.5$, $\varepsilon=0.005$. The incident angle $\theta=3\pi/8$.
Bottom: zoomed view of the two Fano-type transmission anomalies  near $k=3.12$ and $6.24$.
Both Fano resonances occur in the region $D_2$, when there are two and three diffraction orders respectively.
}\label{fig:trans_fano2}
\end{center}
\end{figure}

\begin{figure}
\begin{center}
\hspace*{-2cm}
\includegraphics[height=4.8cm,width=18cm]{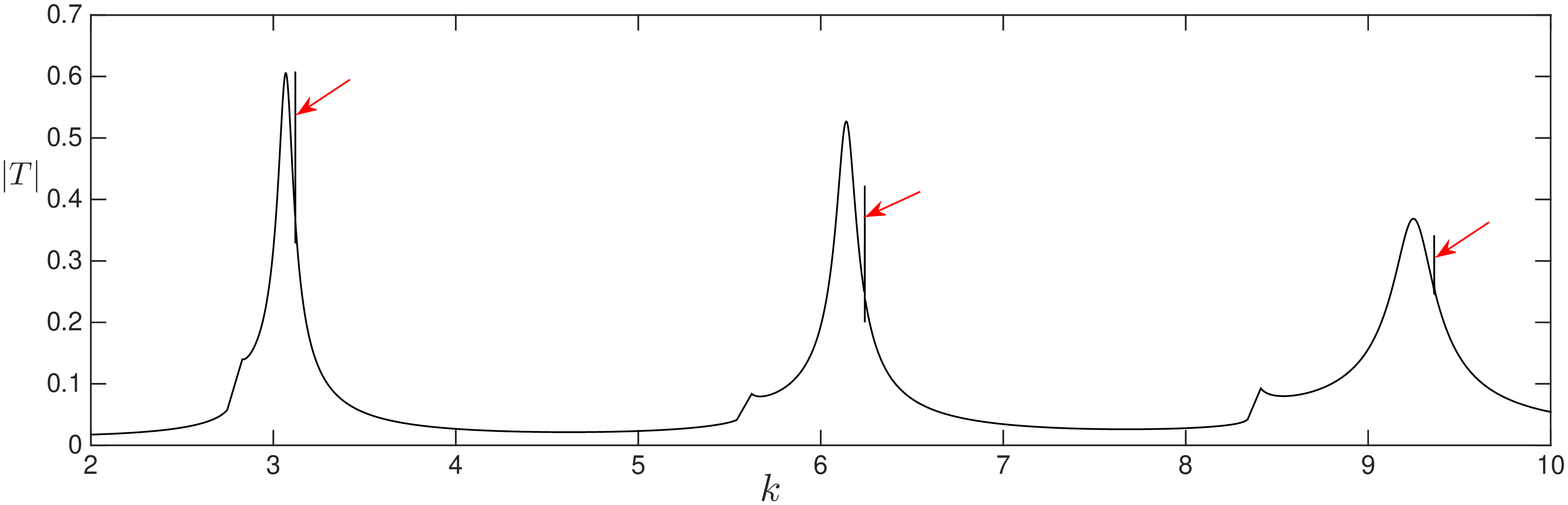}
\hspace*{-1cm}
\includegraphics[height=4.5cm]{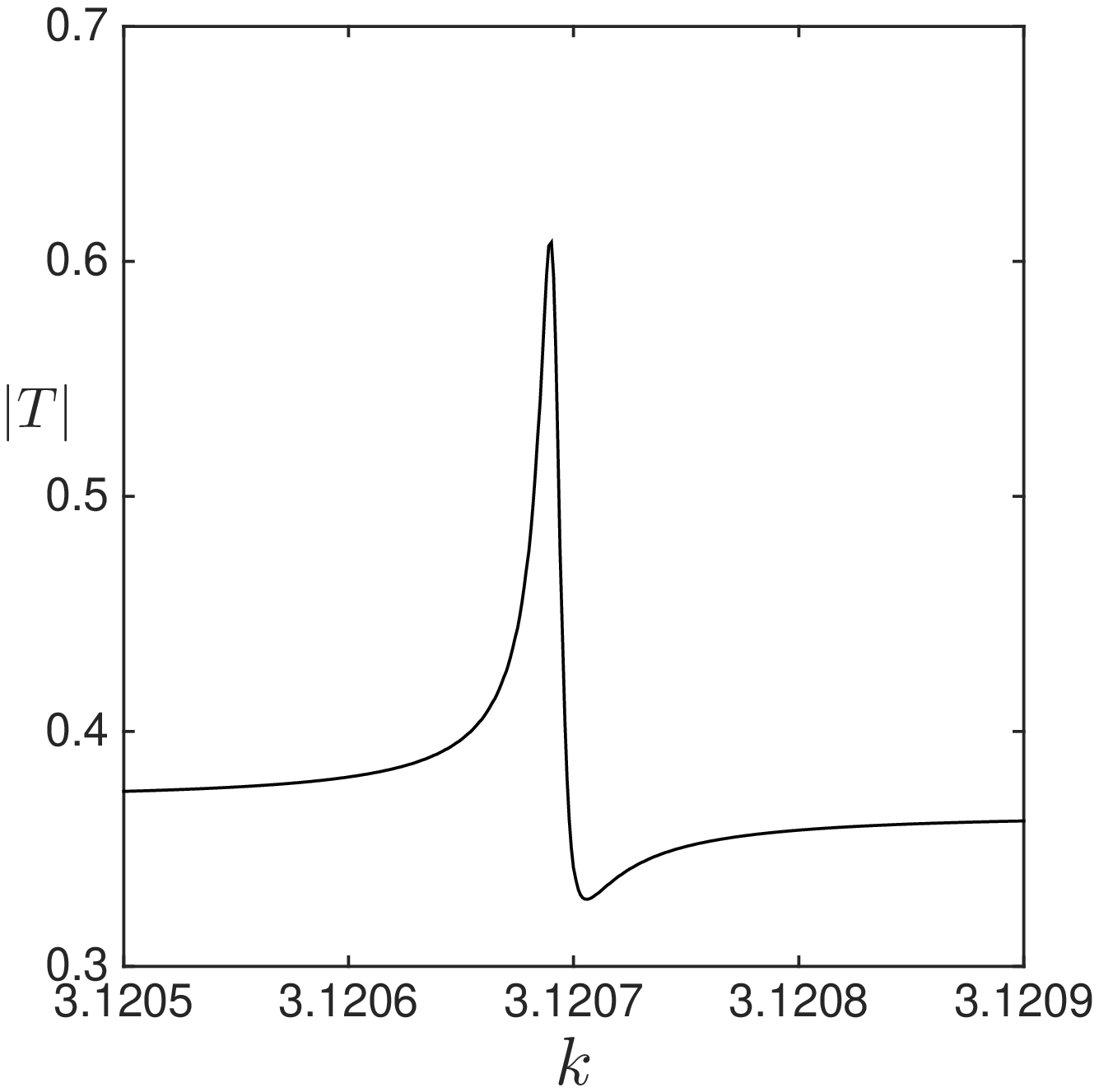}  \quad\quad
\includegraphics[height=4.5cm]{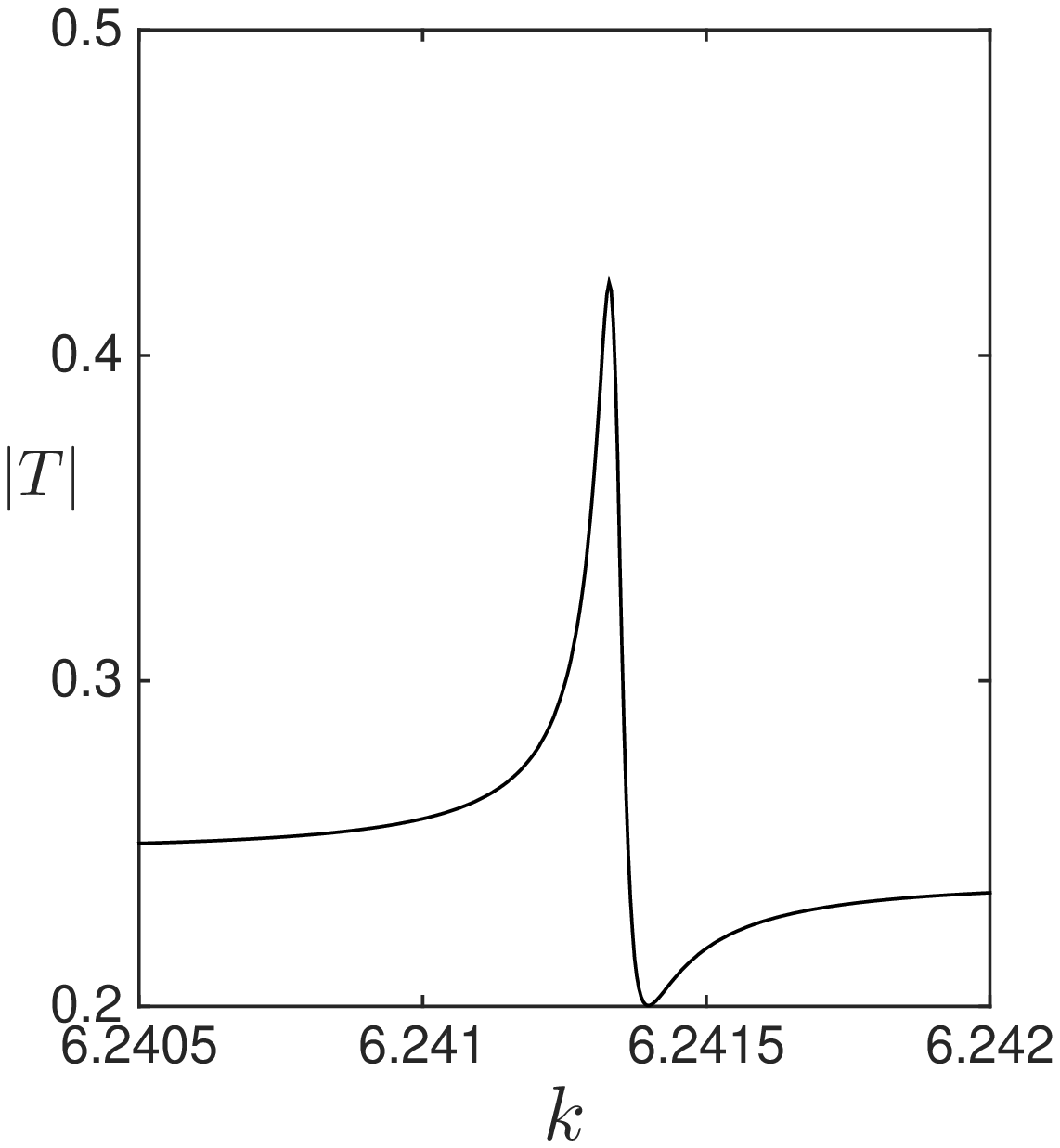}   \quad\quad
\includegraphics[height=4.5cm]{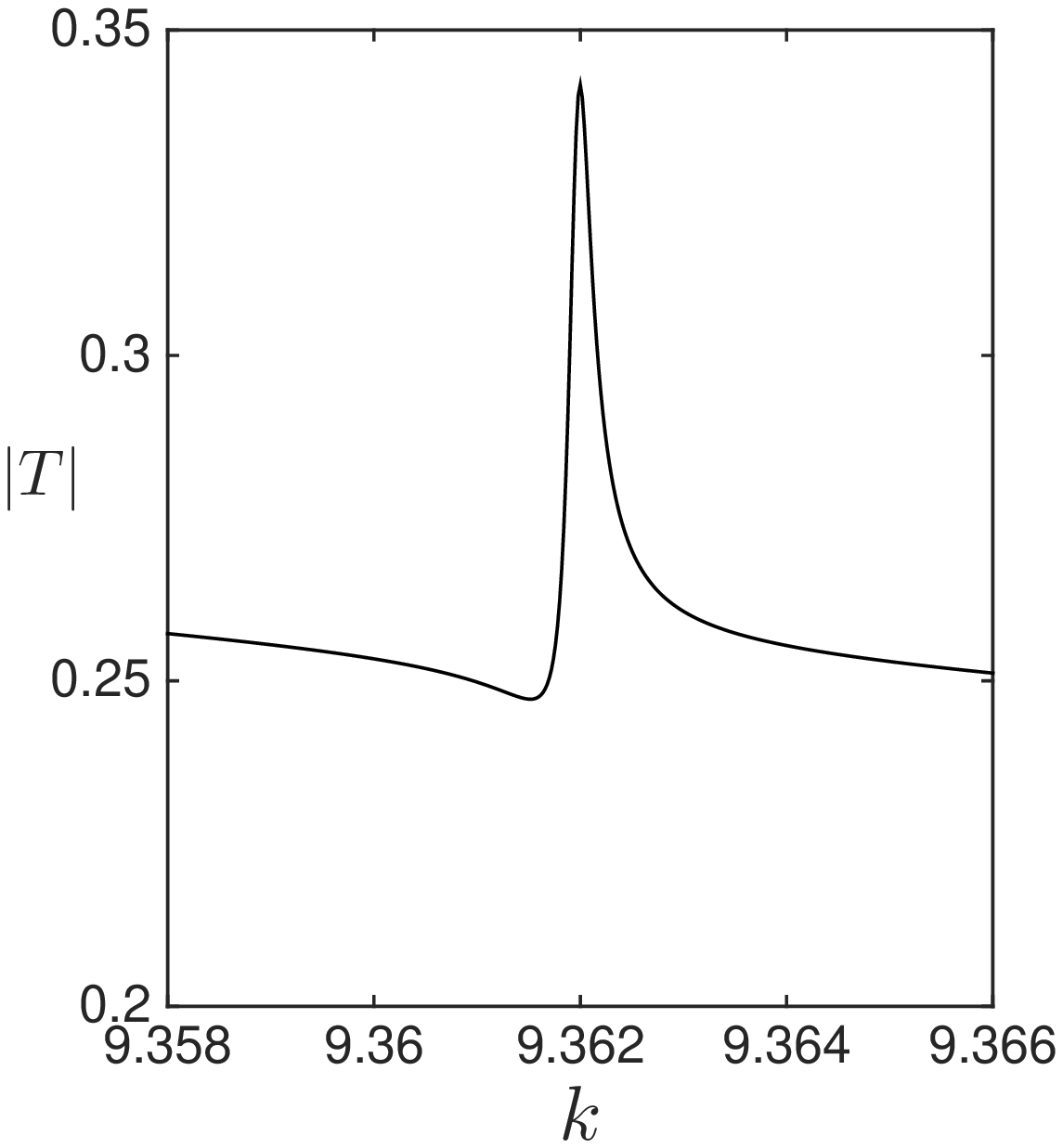}
\caption{Top: Transmission $|T|$ for $k \in (2,10)$ when $d=1.5$, $\varepsilon=0.005$. The incident angle $\theta=\pi/6$.  
Bottom: zoomed view of the three Fano-type transmission anomalies  near $k=3.12$, $6.24$ and $9.36$.
All the three Fano resonances occur in the region $D_2$, when there are two, three and five diffraction orders respectively.
}\label{fig:trans_fano3}
\end{center}
\end{figure}

From the above discussion and Remark \ref{rmk:field}, if 
one substitutes \eqref{t_hat} into the formula \eqref{eq:tn} and use the relation \eqref{eq:sum_w1w2}, then 
for each diffraction order, the transmission coefficient $t_n$ can be viewed as a sum of two signals: one is generated from 
a periodic array of monopoles  (the equivalent charges over the two slit apertures attain the same sign) and the other from a periodic array of dipoles 
(the equivalent charges over the two slit apertures attain opposite signs).
The former almost remains constant near $k_*$, while the phase of the latter changes rapidly. 
Such a combination leads to the Fano-type anomaly in the transmission magnitude $|T|$
as shown in Figure \ref{fig:trans_fano1} - \ref{fig:trans_fano3} for several configurations.

\begin{rmk}\label{rmk:fabry-perot}
In the vicinity of the resonance frequency $Re \, k_m^{(1)}$,  there holds $\varepsilon \hat r_\pm = -\dfrac{\varepsilon \eta_{+}}{\lambda_{1,+}} + O(\varepsilon)$
and $\varepsilon \hat t_\pm = -\dfrac{\varepsilon \eta_{+}}{\lambda_{1,+}} + O(\varepsilon)$. In this scenario, the transmission from the monopole resonant mode is dominant. 
This is the so called Fabry-Perot resonance, which would yield a Lorentzian spectral line shape as shown in Figure \ref{fig:trans_fano1} - \ref{fig:trans_fano3}.
\end{rmk}

Finally, we present several numerical examples. 
We first set the period $d=1.3$ and the slit aperture size $\varepsilon=0.02$. Let the incident angle be $\theta=\pi/6$.
Figure \ref{fig:trans_fano1} demonstrates the leading-order of transmission $|T|$ in \eqref{eq:RT} for $k \in (2,7)$,
where the transmission coefficient $\hat t_\pm := \langle \varphi_2^\pm, 1 \rangle $ is obtained by solving the linear system \eqref{eq:linear_system}.
Fano resonance occurs near $k=3.06$ and $6.12$, which is located in $D_1$ and $D_2$ respectively.
Note that the Lorentzian line shape near $k=2.8$ and $5.9$ is the Fabry-Perot type resonance,
which, in contrast to the Fano resonance involving inference of two resonances, is intrigued by resonance $k_m^{(1)}$ only.
On the other hand, the non-smooth transmission line around $k=3.2$ and $6.4$ corresponds the Rayleigh anomaly, when certain evanescent diffraction modes
convert into propagation modes.

We next set the period $d=1.5$ and the slit aperture size $\varepsilon=0.005$, the transmission for the incident angle $\theta=3\pi/8$ for is shown in Figure \ref{fig:trans_fano2}.
Two Fano resonances appear near $k=3.12$ and $6.24$, and both are in the higher diffraction continuum $D_2$.
Figure \ref{fig:trans_fano3}  plots the transmission for $k \in (2,10)$ when the incident angle becomes smaller.
As predicted by the theory, three Fano resonances occur in this frequency band,
which are all located in the region $D_2$.

\appendix \section{}
\subsection{Asymptotic expansion of $G^{\mathrm{e}}_\varepsilon(X,Y)$}
% Let $\mathbb{N}$ be the set of all nonzero integers,
% $$\mathbb{N}_1 := \{ n \in \mathbb{N} \,;\,  0 < \kappa_n^2 - k^2 < 2|bn|^2  \} 
%   \quad \mbox{and} \quad \mathbb{N}_2 := \mathbb{N} \backslash \mathbb{N}_1.  $$
% It is clear that $\mathbb{N}_1$ is a infinite set, while $\mathbb{N}_2$ is a fintie set.

Let $\tilde{\mathbb{Z}}_1 := \mathbb{Z}_1 \backslash \{ 0 \}$ and $\tilde{\mathbb{Z}}_2 := \mathbb{Z}_2 \backslash \{ 0 \}$,
where $\mathbb{Z}_1$ and $\mathbb{Z}_2$ are defined in \eqref{eq:Z1Z2}.
We split the Green's function as three parts by letting
\begin{equation}\label{eq:Ge_split}
  G^\mathrm{e}_\varepsilon(X,Y) \;=\;  \frac{e^{i\kappa\varepsilon Z}}{d} \left(\frac{1}{i\zeta_0(\kappa,k)}  +  I_1 + I_2\right),
\end{equation}
where $Z=X-Y$, and
\begin{equation}\label{eq:I1_I2}
  I_1:=  \sum_{n\in\tilde{\mathbb{Z}}_1} \frac{e^{i n b\varepsilon Z}}{i\zeta_n(\kappa,k)}, \quad
  I_2:=  \sum_{n\in\tilde{\mathbb{Z}}_2} \frac{e^{i n b\varepsilon Z}}{i\zeta_n(\kappa,k)}. 
\end{equation}

For $n\in\tilde{\mathbb{Z}}_1$, applying the Taylor expansion gives
\begin{equation}\label{eq:I1}
 I_1 =  \sum_{n\in\tilde{\mathbb{Z}}_1} \left( -\frac{1}{|bn|} + \frac{\mathrm{sgn}(n) \kappa}{|bn|^2} \right)  e^{i n b\varepsilon Z} 
+  \sum_{n\in\tilde{\mathbb{Z}}_1} \left(\frac{1}{i \zeta_n(\kappa,k)}  + \frac{1}{|bn|} - \frac{\mathrm{sgn}(n)\kappa}{|bn|^2} \right) + O(\varepsilon).
\end{equation}
For $n\in\tilde{\mathbb{Z}}_2$,  one can write $i\zeta_n = -\sqrt{(\kappa+bn)^2-k^2\,}$, and it follows that
\begin{equation}
  \frac{1}{i\zeta_n(\kappa,k)} \;=\; -\frac{1}{a|n|}
  \left(\sqrt{1+\frac{2\kappa}{bn}+\frac{\kappa^2-k^2}{(bn)^2}\,}\right)^{-1}.
\end{equation}
There holds
\begin{equation*}
  \left| \frac{2\kappa}{bn}+\frac{\kappa^2-k^2}{(bn)^2} \right|  < 1
  \quad \text{for} \; n\in\tilde{\mathbb{Z}}_2,
\end{equation*}
and the Taylor expansion yields
\begin{equation}\label{eq:I2_Taylor}
I_2 =  -\sum_{n\in\tilde{\mathbb{Z}}_2} \frac{1}{b|n|}  
   \left( 1+\sum_{m=1}^\infty c_m \left(\frac{2\kappa}{bn}+\frac{\kappa^2-k^2}{(bn)^2}\right)^m \right)e^{i n b\varepsilon Z},
\end{equation}
where $c_m=-\dfrac{1\cdot3\cdots (2m-1)}{2^m m!}$. Thus we may rewrite $I_2$ as
\begin{equation}\label{eq:I2}
 I_2 = \sum_{n\in\tilde{\mathbb{Z}}_2} \left(-\frac{1}{|bn|} + \frac{\mathrm{sgn}(n)\kappa}{|bn|^2} \right)   e^{i n b\varepsilon Z} - 
 \sum_{n\in\tilde{\mathbb{Z}}_2}  \frac{1}{|bn|}  \sum_{m=2}^\infty \frac{h_m(\kappa,k)}{(bn)^m} e^{i n b\varepsilon Z} 
 \end{equation}
for certain functions $h_m(\kappa,k)$ satisfying
\begin{equation}\label{eq:sum_series1}
 -\frac{1}{|bn|} \sum_{m=2}^\infty \frac{h_m(\kappa,k)}{(bn)^m} = \frac{1}{i\zeta_n(\kappa,k)} + \frac{1}{|bn|} - \frac{\mathrm{sgn}(n)\kappa}{|bn|^2}, \quad n\in\mathbb{Z}_2. 
\end{equation}

\begin{itemize}
\item [(i)] $(\kappa,k) \in D_1$.  It follows that $\tilde{\mathbb{Z}}_1=\emptyset$,
$\tilde{\mathbb{Z}}_2=\mathbb{Z}\backslash\{0\}$, and $I_1$ vanishes.
If $\kappa =0$, an examination of the expansion \eqref{eq:I2_Taylor} reveals that $h_m(k,0)=0$ when $m$ is odd. The following relations hold
(cf.~\cite{collin, kress})
\begin{eqnarray}
 - \sum_{n \neq 0 }\dfrac{1}{|n|} e^{inb\varepsilon Z} &=& \ln\left(4\sin^2\frac{b\varepsilon Z}{2}\right),  \label{eq:sum_eq1} \\
\sum_{n \neq 0 }^\infty \dfrac{1}{|n| n^m} e^{ina\varepsilon Z} &=& \sum_{n \neq 0}^\infty \dfrac{1}{|n| n^m}+ i^m \cdot O\left(\varepsilon^{m}Z^m\ln(|\varepsilon Z|)\right)  \quad (m \ge 2), \label{eq:sum_eq2}
\end{eqnarray}
where the $O\left(\varepsilon^{m}Z^m\ln(|\varepsilon Z|)\right)$ is real. By substituting \eqref{eq:sum_eq1} and \eqref{eq:sum_eq2} with even $m$ into \eqref{eq:I2}, 
and using the relation \eqref{eq:sum_series1}, we obtain
\begin{eqnarray}
I_2 &=&  \frac{1}{b}\ln\left(4\sin^2\frac{b\varepsilon Z}{2}\right)  -\sum_{m=2}^\infty \frac{h_m(\kappa,k)}{b^{m+1}} \sum_{n\neq 0}  \frac{1}{|n| n^m} + O(\varepsilon^2|\ln\varepsilon|) \nonumber  \\
  &=&   \frac{2}{b} \ln (b \varepsilon |Z|)  + \sum_{n\neq 0}  \left( \frac{1}{i\zeta_n(k,0)} + \frac{1}{|bn|}  \right) + O(k^2\varepsilon^2|\ln\varepsilon|).   \label{eq:I2_exp1}
\end{eqnarray} 
The desired asymptotic expansion follows by substituting $I_2$ into \eqref{eq:Ge_split}.
In addition, the above expansion shows that $r_\mathrm{e}(0,\varepsilon;X, Y)=O(\varepsilon^2|\ln\varepsilon|)$ is a function of $|Z|:=|X-Y|$ and is real when $\kappa=0$.

If $\kappa\neq 0$,  from the relation \eqref{eq:sum_eq2},
the second term in \eqref{eq:I2} can be recast as
$$ -\sum_{m=2}^\infty \frac{h_m(\kappa,k)}{b^{m+1}} \sum_{n\neq 0}  \frac{1}{|n| n^m} + O(\varepsilon^2|\ln\varepsilon|) + O(\kappa\varepsilon^3|\ln\varepsilon|),  $$
where $O(\varepsilon^2|\ln\varepsilon|)$ and $O(\kappa\varepsilon^3|\ln\varepsilon|)$ represent the high-order terms arising from the sum of the even and odd $m$ respectively.
Using \eqref{eq:sum_series1} again, it follows that
\begin{eqnarray*}
I_2 &=& \sum_{n \neq 0 } \left( -\frac{1}{|bn|} + \frac{\mathrm{sgn}(n) \kappa }{|bn|^2} \right)  e^{i n b\varepsilon Z} + 
\sum_{n \neq 0} \left(\frac{1}{i \zeta_n(\kappa,k)}  + \frac{1}{|bn|} - \frac{\mathrm{sgn}(n)\kappa}{|bn|^2} \right) \\
&& + O(\varepsilon^2|\ln\varepsilon|) + O(\kappa\varepsilon^3|\ln\varepsilon|).
\end{eqnarray*}
Now an application of the relation
\begin{equation*}
  \sum_{n \neq 0 }^\infty \frac{\mathrm{sgn}(n)}{|n|^{2}} e^{inb\varepsilon Z} = - i (2b \varepsilon Z \ln (b \varepsilon |Z|) + O(\varepsilon).
\end{equation*}
yields
\begin{equation}\label{eq:I2_exp2}
I_2 =  (1 - i \kappa\varepsilon Z)  \frac{2}{b} \ln (b \varepsilon |Z|) + \sum_{n \neq 0} \left(\frac{1}{i \zeta_n(\kappa,k)}  + \frac{1}{|bn|} \right) + O(\varepsilon^2|\ln\varepsilon|) +O(\kappa\varepsilon)+O(\kappa\varepsilon^3|\ln\varepsilon|).  
\end{equation}
By substituting the above into \eqref{eq:Ge_split} and using the Taylor expansion for $e^{i\kappa\varepsilon Z}$, we obtain the desired asymptotic expansion for $G_e(X,Y)$. In addition, by examing the above calculations, 
it is clear that $r_\mathrm{e}(\kappa,\varepsilon; X,Y) = r_\mathrm{e}(0,\varepsilon; X,Y) + O(\kappa\varepsilon) $,
since the $O(\varepsilon^2|\ln\varepsilon|)$ term in \eqref{eq:I2_exp1} and \eqref{eq:I2_exp2} represents the high-order term 
arising from the sum of the even $m$ and are the same for $\kappa=0$ and $\kappa \neq 0$.

\item [(ii)] $(\kappa,k) \in D_2$. For each $m \ge 2$, the Taylor expansion gives
$$ \sum_{n \in \tilde{\mathbb{Z}}_1} \dfrac{1}{|n| n^m} e^{inb\varepsilon Z} = \sum_{n\in \tilde{\mathbb{Z}}_1}^\infty \dfrac{1}{|n| n^m}+  O(\varepsilon). $$
Thus by applying \eqref{eq:sum_eq2}, we see that
$$  \sum_{n\in\tilde{\mathbb{Z}}_2}  \frac{1}{|n| n^m}e^{i n b\varepsilon Z} = \sum_{n\neq0}  \frac{1}{|n| n^m}e^{i n b\varepsilon Z} -
\sum_{n\in\tilde{\mathbb{Z}}_1}  \frac{1}{|n| n^m}e^{i n b\varepsilon Z}  =  \sum_{n\in\tilde{\mathbb{Z}}_2}  \frac{1}{|n| n^m} + O(\varepsilon^2|\ln\varepsilon|) +O(\varepsilon). $$
As such the second term in \eqref{eq:I2} can be expressed as
$$  -\sum_{m=2}^\infty \frac{h_m(\kappa,k)}{b^{m+1}} \sum_{n\in\tilde{\mathbb{Z}}_2}  \frac{1}{|n| n^m} + O(\varepsilon^2|\ln\varepsilon|) +  O(\varepsilon),  $$
and we obtain
\begin{equation}\label{eq:I2_exp3}
I_2 = \sum_{n\in\tilde{\mathbb{Z}}_2} \left( -\frac{1}{|bn|} + \frac{\mathrm{sgn}(n) \kappa }{|bn|^2} \right)  e^{i n b\varepsilon Z} + 
\sum_{n\in\tilde{\mathbb{Z}}_2} \left(\frac{1}{i \zeta_n(\kappa,k)}  + \frac{1}{|bn|} - \frac{\mathrm{sgn}(n)\kappa}{|bn|^2} \right) + O(\varepsilon).
\end{equation}
A combination of \eqref{eq:I1} and \eqref{eq:I2_exp3} leads to
$$
I_1 + I_2 = \sum_{n \neq 0 } \left( -\frac{1}{|bn|} + \frac{\mathrm{sgn}(n) \kappa }{|bn|^2} \right)  e^{i n b\varepsilon Z} + 
\sum_{n \neq 0} \left(\frac{1}{i \zeta_n(\kappa,k)}  + \frac{1}{|bn|} - \frac{\mathrm{sgn}(n)\kappa}{|bn|^2} \right) + O(\varepsilon).
$$
The desired asymptotic expansion follows by the same calculations as in case (i).
\end{itemize}

\bibliography{references}

\begin{thebibliography}{99}
\bibitem{habib-book}
{\sc H. Ammari, H. Kang, B. Fitzpatrick, M. Ruiz, S. Yu and H. Zhang}
{\em Mathematical and Computational Methods in Photonics and Phononics}, 
Mathematical Surveys and Monographs, Volume 235, American Mathematical Society, Providence, 2018.


\bibitem{habib2019}
{\sc H. Ammari, B. Fitzpatrick, H. Lee, S. Yu and H. Zhang}
{\em Double-negative acoustic metamaterials},
Quarterly of Applied Mathematics, \textbf{77} (2019), no. 4, 767-791.


\bibitem{bao95}
{\sc G. Bao, D. Dobson, and Cox}, {\em Mathematical studies in rigorous grating theory}, J. Opt. Soc. Amer. A, \textbf{12} (1995), 1029-1042.

\bibitem{bonnet_starling94}
{\sc  A. Bonnet-Bendhia and F. Starling}, {\em Guided waves by electromagnetic gratings and non-uniqueness examples for the diffraction problem},
Math. Meth. Appl. Sci., \textbf{17} (1994), 305-338.

\bibitem{eric10-2}
{\sc  J. F. Babadjian, E. Bonnetier and F. Triki},
{\em Enhancement of electromagnetic fields caused by interacting subwavelength cavities},
Multiscale Model. Simul., \textbf{8} (2010), 1383-1418. 

\bibitem{bugakov}
{\sc E. Bulgakov and Almas F. Sadreev}, {\em Bloch bound states in the radiation continuum in a periodic array of dielectric rods},
Phy. Re. A, \textbf{90} (2014): 053801.

\bibitem{collin}
{\sc R. Collin}, Field Theory of Guided Waves (2nd Edition), Wiley-IEEE Press, 1990.

\bibitem{fano}
{\sc U. Fano}, {\em Effects of configuration interaction on intensities and phase shifts}, Phy. Rev., \textbf{124} (1961), 1866.

\bibitem{fan1}
{\sc S. Fan}, {\em Sharp asymmetric line shapes in side-coupled waveguide-cavity systems}, Appl. Phy. Lett., \textbf{80} (2002): 908.

\bibitem{fan2}
{\sc S. Fan, W. Suh, and J. D. Joannopoulos},
{\em Temporal coupled-mode theory for the {F}ano resonance in optical resonators},
J. Opt. Soc. Am. A, \textbf{20}(3) (2003) 569--572.


\bibitem{hsu1}
{\sc C. W. Hsu, B. Zhen, S. L. Chua, S. Johnson, J. Joannopoulos, and M. Soljacic},
{\em Bloch surface eigenstates within the radiation continuum},
Light: Science App., 2(e84 doi:10:1038/Isa.2013.40), 2013.

\bibitem{hsu2}
{\sc C. W. Hsu, B. Zhen, J. Lee, S. L. Chua, S. Johnson, J. Joannopoulos, and M. Soljacic},
{\em Observation of trapped light within the radiation continuum},
Nature, (doi:10.1038/nature12289), July 2013.

\bibitem{hsu3}
{\sc C. W. Hsu, \textit{et al.}}, {\em Bound states in the continuum}, Nat. Rev. Mater., \textbf{1} (2016), 16048.

\bibitem{kress}
{\sc R. Kress}, {\em Linear Integral Equations}, Applied Mathematical Sciences, vol. 82, Springer-Verlag, Berlin, 1999.

\bibitem{li}
B. Li, \textit{et al.}, {\em Experimental controlling of Fano resonance in indirectly coupled whispering-gallery microresonators}, Appl. Phy. Lett., \textbf{100} (2012), 021108.

\bibitem{lin_zhang17}
{\sc  J. Lin and H. Zhang}, {\em Scattering and field enhancement of a perfect conducting narrow slit}, SIAM J. Appl. Math.,  (2017), 951--976.

\bibitem{lin_zhang18_1}
{\sc  J. Lin and H. Zhang}, {\em Scattering by a periodic array of subwavelength slits I: field enhancement in the diffraction regime}, Multiscale Model. Simul., \textbf{16} (2018), 922--953.


\bibitem{lin_shipman_zhang}
{\sc  J. Lin, S. Shipman, and H. Zhang}, {\em A mathematical theory for Fano resonance in a periodic array of narrow slits}, SIAM J. Appl. Math., to appear.


\bibitem{limonov}
{\sc M. Limonov, \textit{et al.}},  {\em Fano resonances in photonics}, Nat. Photon. \textbf{11} (2017): 543.
 
\bibitem{linton98}
{\sc C. Linton}, {\em the Green function for the two-dimensional Helmholtz equation in periodic domains}, J. Eng. Math. \textbf{33} (1998), 377--401.

\bibitem{lukyanchuk}
{\sc B. Luk'yanchuk, \textit {et al.}}, {\em The Fano resonance in plasmonic nanostructures and metamaterials}, Nat. Mat. \textbf{9} (2010): 707.

\bibitem{porter}
{\sc R. Porter and D. V. Evans}, {\em Embedded Rayleigh-Bloch surface waves along periodic rectangular arrays}, Wave Motion \textbf{43} (2005): 29-50.

\bibitem{rayleigh}
L. Rayleigh, {\em On the dynamical theory of gratings}, P. R. Soc. London,  \textbf{79} (1907): 399-416.

\bibitem{shipman05}
{\sc S. P. Shipman and S. Venakides},  {\em Resonant transmission near nonrobust periodic slab modes}, Phy. Rev. E, \textbf{71} (2005), 026611.

\bibitem{shipman10}
{\sc S. P. Shipman}, {\em Resonant scattering by open periodic waveguides},  Chapter 2 in Wave Propagation in Periodic Media: Analysis, Numerical Techniques and Practical Applications, M. Ehrhardt, ed., E-Book Series PiCP, Bentham Science Publishers, Vol. 1 (2010).

\bibitem{TKT}
{\sc K. Totsuka, N. Kobayashi, and M. Tomita}, {\em Slow light in coupled-resonator-induced transparency}, Phy. Rev. Lett., \textbf{98} ( 2007), 213904.

 \bibitem{lu1}
{\sc L. Yuan and Y-Y. Lu}, {\em Propagating Bloch modes above the lightline on a periodic array of cylinders},  J. Phy. B: Atomic, Molecular and Optical Physics, \textbf{50} (2017), 05LT01.

\bibitem{lu2}
{\sc L. Yuan and Y-Y. Lu}, {\em Bound states in the continuum on periodic structures: perturbation theory and robustness}, Opt. Lett., \textbf{42} (2017), 4490--4493.

\bibitem{zhou}
L. Zhou and A. W. Poon, {\em Fano resonance-based electrically reconfigurable add-drop filters in silicon microring resonator-coupled Mach-Zehnder interferometers},
Opt. Lett., \textbf{32} (2007), 781-783.

\end{thebibliography}

\end{document}